\documentclass{article}
\usepackage[utf8]{inputenc}
\usepackage{amsmath,amsfonts,amsthm,geometry,import,mathrsfs}
\usepackage{xcolor}
\RequirePackage[colorlinks,linkcolor=blue,citecolor=blue,urlcolor=blue]{hyperref}
\usepackage[english]{babel}
\usepackage{comment}

\newcommand{\R}{\mathbb{R}}

\newcommand{\RD}{\mathbb{R}^d}
\newcommand{\N}{\mathbb{N}}

\newcommand{\cA}{\mathcal{A}}

\newcommand{\cC}{\mathcal{C}}
\newcommand{\cF}{\mathcal{F}}

\newcommand{\cH}{\mathcal{H}}
\newcommand{\cI}{\mathcal{I}}

\newcommand{\cP}{\mathcal{P}}

\newcommand{\cT}{\mathcal{T}}

\newcommand{\cV}{\mathcal{V}}
\newcommand{\cW}{\mathcal{W}}

\newcommand{\scrI}{\mathscr{I}}

\newcommand{\bes}{\begin{equation*}}
\newcommand{\ees}{\end{equation*}}
\newcommand{\beas}{\begin{eqnarray*}}
\newcommand{\eeas}{\end{eqnarray*}}
\newcommand{\bea}{\begin{eqnarray}}
\newcommand{\eea}{\end{eqnarray}}
\newcommand{\be}{\begin{equation}}
\newcommand{\ee}{\end{equation}}
\newcommand{\bei}{\begin{itemize}}
\newcommand{\eei}{\end{itemize}}
\newcommand{\bec}{\begin{cases}}
\newcommand{\eec}{\end{cases}}
\newcommand{\ben}{\begin{enumerate}}
\newcommand{\een}{\end{enumerate}}

\newcommand{\bbE}{\mathbb{E}}

\newcommand{\bbl}{\begin{block}}
\newcommand{\ebl}{\end{block}}

\newcommand{\De}{\mathrm{d}}
\newcommand{\curv}{\kappa}

\newcommand{\rmH}{\mathrm{H}}
\newcommand{\rmP}{\mathrm{P}}
\newcommand{\rmQ}{\mathrm{Q}}
\newcommand{\rmR}{\mathrm{R}}

\newtheorem{mydef}{Definition}[section]
\newtheorem{prop}{Proposition}[section]
\newtheorem{theorem}{Theorem}[section]
\newtheorem{lemma}{Lemma}[section]

\newtheorem{remark}{Remark}[section]

\newtheorem{cor}{Corollary}[section]
\newcommand{\absnent}{\tilde{\mathcal{F}}}
\newcommand{\nent}{\mathcal{F}}
\newcommand{\tang}{\mathrm{Tan}}
\newcommand{\nonlinentcost}{\mathscr{C}_T}
\newcommand{\cq}{\mathscr{E}_{\rmP}}
\newcommand{\tansp}{\cT_{\mu}\cP_2}
\newcommand{\tanspt}{\cT_{\mu_t}\cP_2}
\newcommand{\otgrad}{\mathrm{grad}^{\mathcal{W}_2}}
\newcommand{\othess}{\mathrm{Hess}^{\mathcal{W}_2}}
\newcommand{\mui}{\mu^{\mathrm{in}}}
\newcommand{\muf}{\mu^{\mathrm{fin}}}
\newcommand{\nonlinfish}{\mathcal{I}_{\nent}}
\newcommand{\leb}{\lambda}
\newcommand{\mkv}{\rmP^{\text{\tiny{MKV}}}}

\title{The mean field Schr\"odinger problem: ergodic behavior, entropy estimates and functional inequalities}
\author{Julio Backhoff\thanks{Faculty of Mathematics, University of Vienna, Oskar-Morgenstern-Platz 1, 1090 Vienna, Austria. julio.backhoff@univie.ac.at}, Giovanni Conforti\thanks{D\'epartement de Math\'ematiques Appliqu\'ees, \'Ecole Polytechnique, Route de Saclay, 91128, Palaiseau Cedex, France. giovanni.conforti@polytechnique.edu }, Ivan Gentil\thanks{Institut Camille Jordan, Univ Claude Bernard Lyon 1,  France. gentil@math.univ-lyon1.fr}, Christian L\'eonard\thanks{Modal’X, UPL, Univ Paris Nanterre, F92000 Nanterre France. christian.leonard@u-paris10.fr}}

\begin{document}

\maketitle

\begin{abstract} 
{We study the mean field Schr\"odinger problem (MFSP), that is the problem of finding the most likely evolution of a cloud of \emph{interacting} Brownian particles conditionally on the observation of their initial and final configuration. Its rigorous formulation is in terms of an optimization problem with marginal constraints whose objective function is the large deviation rate function associated with a system of weakly dependent Brownian particles. We undertake a fine study of the dynamics of its solutions,  including quantitative energy dissipation estimates yielding the exponential convergence to equilibrium as the the time between observations grows larger and larger, as well as a novel class of functional inequalities involving the mean field entropic cost (i.e. the optimal value in (MFSP)). Our strategy unveils an interesting connection between forward backward stochastic differential equations and the Riemannian calculus on the space of probability measures introduced by Otto, which is of independent interest.}
\end{abstract}

\tableofcontents

\section{Introduction and statement of the main results}

In the seminal works \cite{Schr,Schr32} E.\,Schr\"odinger addressed the problem of finding the most likely evolution of a cloud of independent Brownian particles conditionally on the observation of their initial and final configuration. In modern language this is an entropy minimization problem with marginal constraints.
The aim of this work is to take the first steps in the understanding of the Mean Field Schr\"odinger Problem, obtained by replacing in the above description the \emph{independent} particles by \emph{interacting} ones. 

To obtain an informal description of the problem, consider $N$ Brownian particles $(X^{i,N}_t)_{t\in[0,T],1\leq i \leq N}$ interacting through a pair potential $W$

\be\label{partsys intro}
\left\{ 
\begin{array}{rl}
     \De X_t^{i,N}=&-\frac{1}{N} \sum_{k=1}^N\nabla W(X_t^{i,N}-X_t^{k,N})\De t+\De B^i_t  \\
     X_0^{i,N}\sim& \mu^{\mathrm{in}}.
\end{array}
\right .
\ee
Their evolution is encoded in the random empirical path measure 
\be\label{random empirical path}  \textstyle  \frac{1}{N}\sum_{i=1}^N\delta_{X^{i,N}_{\cdot}} .\ee
At a given time $T$, the configuration of the particle system is visible to an external observer that finds it close to an ``unexpected" (\emph{\'ecart spontan\'e et considerable} in \cite{Schr32}) probability measure $\muf$, namely
\be\label{Schr obs}
\frac{1}{N}\sum_{i=1}^N\delta_{X^{i,N}_{T}} \approx \muf
\ee
It is a classical result \cite{Ta84,DaGae87,BuDuFi12} that the sequence of empirical path measures \eqref{random empirical path} obeys the large deviations principle (LDP). Thus, the problem of finding the most likely evolution conditionally on the observations is recast as the problem of minimizing the LDP rate function among all path measures whose marginal at time $0$ is $\mui$ and whose marginal at time $T$ is $\muf$. This is the mean field Schr\"odinger problem (MFSP). Extending naturally the classical terminology we say that an optimal path measure is a mean field Schr\"odinger bridge (henceforth MFSB) and the optimal value is the mean field entropic cost. The latter generalizes both the Wasserstein distance and the entropic cost. 

The classical Schr\"odinger problem has been the object of recent intense research activity (see \cite{LeoSch}). This is due to the computational advantages deriving from introducing an entropic penalization in the Monge-Kantrovich problem \cite{cuturi2013sinkhorn} or to its relations with functional inequalities, entropy estimates and the geometrical aspects of optimal transport. Our article contributes to this second line of research, recently explored by the papers  \cite{conforti2018second,gentil2018dynamical,leger2019hopf,ripani2017convexity,karatzas2018pathwise}. Leaving all precise statements to the main body of the introduction, let us give a concise summary of our contributions.

\paragraph{Dynamics of mean field Schr\"odinger bridges}  Our mean field version of the Schr\"odinger problem stems from fundamental results in large deviations for weakly interacting particle systems such as \cite{Ta84,DaGae87} and shares some analogies with the control problems considered in \cite{chen2018steering} and with the article \cite{benamou2018entropy} in which an entropic formulation of second order variational mean field games is studied.  Among the more fundamental results we establish for the mean field Schrödinger problem, we highlight
\begin{itemize}
    \item the existence of MFSBs and, starting from the original large deviations formulation, the derivation of both an equivalent reformulation in terms of a McKean-Vlasov control problem as well as a Benamou-Brenier formula,
    \item establishing that MFSBs solve forward backward stochastic differential equations (FBSDE) of McKean-Vlasov type (cf.\ \cite{CaDe13,CaDe15}).
\end{itemize}
The proof strategy we adopt in this article combines ideas coming from large deviations and stochastic calculus of variations, see \cite{CrLa15, ustunel2007sufficient,emery1988cherchant}. 
Another interesting consequence of having a large deviations viewpoint is that we can also exhibit some regularity properties of MFSBs, taking advantage of F\"ollmer's results \cite{follmer1986time} on time reversal. Building on \cite{conforti2018second,gentil2018dynamical} we establish a link between FBSDEs and the Riemannian calculus on probability measures introduced by Otto \cite{otto2001geometry} that is of independent interest and underlies our proof strategies. In a nutshell, the seminal article \cite{jordan1998variational} established that the heat equation is the gradient flow of the relative entropy w.r.t.\ the squared Wasserstein distance. Thus, classical first order SDEs yield probabilistic representations for first order ODEs in the Riemannian manifold of optimal transport. Our observation may be seen as the second order counterpart to the results of \cite{jordan1998variational}: indeed we will present an heuristic strongly supporting the fact that Markovian solutions of ``second order" \emph{pathwise} equations (FBSDEs) yield probabilistic representations for second order ODEs in the Riemannian manifold of optimal transport.

\paragraph{Ergodicity of Schr\"odinger bridges and functional inequalities} Consider again \eqref{partsys intro} and assume that $W$ is convex so that the particle system is rapidly mixing and there is a well defined notion of \emph{equilibrium configuration} $\mu_{\infty}$. If $N$ and $T$ are large, one expects that
\bei
\item[(i)] The configurations $\frac{1}{N}\sum_{i=1}^N\delta_{X^i_t}$ at times $t=0,T/2,T$ are almost independent.
\item[(ii)] The configuration at $T/2$ is with high probability very similar to  $\mu_{\infty}$.
\eei

Because of (i), even when the external observer acquires the information \eqref{Schr obs}, he/she still expects (ii) to hold.
 Thus mean field Schr\"odinger bridges are to spend most of their time around the equilibrium configuration. All our quantitative results originate in an attempt to justify rigorously this claim. 
 
 In this work we obtain a number of precise quantitative energy dissipation estimates. These lead us to the main quantitative results of the article:
 \begin{itemize}
    \item we characterize the long time behavior of MFSBs, proving exponential convergence to equilibrium with sharp exponential rates,
    \item we derive a novel class of functional inequalities involving the mean field entropic cost. Precisely, we obtain a Talagrand inequality and an HWI inequality that generalize those previously obtained in \cite{carrillo2003kinetic} by Carrillo, McCann and Villani.
\end{itemize}
Regarding the second point above, we can in fact retrieve (formally) the inequalities in \cite{carrillo2003kinetic} by looking at asymptotic regimes for the  mean field Schrödinger problem. Besides the intrinsic interest and their usefulness in establishing some of our main results, our functional inequalities may have consequences in terms of concentration of measure and hypercontractivity of non linear semigroups, but this is left to future work.
 
 The fact that optimal curves of a given optimal control problem spend most of their time around an equilibrium is known in the literature as the \emph{turnpike property}. The first turnpike theorems have been established in the 60's for problems arising in econometry \cite{mckenzie1963turnpike}; general results for deterministic finite dimensional problems are by now available, see \cite{trelat2015turnpike}. In view of the McKean-Vlasov formulation of the  mean field Schrödinger problem, some of our 
 results may be viewed as turnpike theorems as well, but for a class of infinite dimensional and stochastic problems. An interesting feature is that, by exploiting the specific structure of our setting, we are able to establish the turnpike property in a quantitative, rather than qualitative form. The McKean-Vlasov formulation also connects our findings with the study of the long time behavior of mean field games \cite{cardaliaguet2012long,cardaliaguet2013long,cardaliaguet2013KAM,cardaliaguet2017long}. 
 
 Concerning the proof methods, our  starting point is Otto calculus and the recent rigorous results of \cite{conforti2018second} together with the heuristics put forward in \cite{gentil2018dynamical}. The first new ingredient of our proof strategy is the above mentioned connection between FBSDEs and Otto calculus that plays a key role in turning the heuristics into rigorous statements. It is worth remarking that using a pathwise approach does not just provide with a way of making some heuristics rigorous, but it also permits to obtain a stronger form of some of the results conjectured in \cite{gentil2018dynamical} which then simply follow by averaging pathwise estimates. { The second new ingredient in our proofs involves a conserved quantity that plays an analogous role to the total energy of a phyisical system. For such quantity we derive a further functional inequality which seems to be novel already in the classical Schrödinger problem (i.e.\ for independent particles) and allows to establish the turnpike property.}

\paragraph{Structure of the article}
In the remainder introductory section we state and comment our main results. In Section \ref{OTFBSDEsection} we provide a geometrical interpretation sketching some interesting heuristic connections between optimal transport and stochastic calculus. The material of this section is not used later on; therefore the reader who is not interested in optimal transport may avoid it. Sections \ref{sec formulations proofs} and \ref{quant res proofs} contain the proofs of our main results, the former being devoted to the results concerning the dynamics of MFSBs and the latter one dealing with the ergodic results. Finally an appendix section contains some technical results.

\subsection{Frequently used notation}
\bei
\item $(\Omega,\mathcal{F}_t,\mathcal{F}_T)$ is the canonical space of $\mathbb R^d$-valued continuous paths on $[0,T]$, so $\{\mathcal{F}_{t}\}_{t\leq T}$ is the coordinate filtration. $\Omega$ is endowed with the uniform topology.
\item $\cP(\Omega)$ and $\cP(\RD)$ denote the set of Borel probability measures on $\Omega$ and $\RD$ respectively.
\item $(X_t)_{t\in[0,T]}$ is the canonical (ie.\ identity) process on $\Omega$.
\item $\rmR^\mu$ is the Wiener measure with starting distribution $\mu$.
\item $\cH(\rmP|\rmQ)$ denotes the relative entropy of $\rmP$ with respect to $\rmQ$, defined as $\bbE_\rmP\left[\log\left(\frac{\De\rmP}{\De\rmQ} \right)\right]$ if $\rmP\ll\rmQ$ and $+\infty$ otherwise. 
\item $\rmP_t$ denotes the marginal distribution of a measure $\rmP\in\cP(\Omega)$ at time $t$.
\item $\cP_{\beta}(\Omega)$ is the set of measures on $\Omega$ for which $\sup_{t\leq T}|X_t|^\beta$ is integrable. $\cP_{\beta}(\RD)$ is the set of measures on $\RD$ for which the function $|\cdot|^\beta$ is integrable.
\item The $\beta$-Wasserstein distance on $\cP_{\beta}(\Omega)$ is defined by
\be \notag\textstyle
\cP_{\beta}(\Omega)^2\ni(\rmP,\rmQ)\mapsto \mathcal  W_\beta(\rmP,\rmQ):=\left( \inf_{Y\sim \rmP,Z\sim \rmQ}\bbE\left[\sup_{t\in[0,T]}|Y_t-Z_t|^\beta \right ] \right)^{1/\beta}.
\ee
With a slight abuse of notation we also denote by $\mathcal W_\beta$ the $\beta$-Wasserstein distance on $\cP_{\beta}(\mathbb R^d)$ defined analogously.
\item  For a given measurable marginal flow $[0,T]\ni t\mapsto \mu_t \in \cP(\RD)$, we denote by $L^2((\mu_t)_{t\in[0,T]})$ the space of square integrable functions from $[0,T]\times\RD$ to $\RD$ associated to the reference measure $\mu_t(\De x)\De t$ and the corresponding almost-sure identification. We consider likewise the Hilbert space $$\rmH_{-1}((\mu_t)_{t\in[0,T]}),$$ defined as the closure in $L^2((\mu_t)_{t\in[0,T]})$
of the smooth subspace 
\bes
 \left\{ \Psi:[0,T]\times \RD \rightarrow \RD \  \text{s.t.} \ \Psi = \nabla \psi,\, \psi \in \cC^{\infty}_c([0,T] \times \RD) \right\}.
 \ees
\item $\gamma$ and $\lambda$ are respectively the standard Gaussian and Lebesgue measure in $\mathbb R^d$.
\item $\cC^{l,m}([0,T]\times\RD;\mathbb R^k)$ is the set of functions from $[0,T]\times\RD $ to $ \mathbb R^k$ which have $l$ continuous derivatives in the first (ie.\ time) variable and $m$ continuous derivatives in the second (ie.\ space) variable. The space $\cC^{m}(\RD;\mathbb R^k)$ is defined in the same way. $\cC^{\infty}_c([0,T]\times\RD)$ is the space of real-valued smooth functions on $[0,T]\times\RD $ with compact support. The gradient $\nabla$ and Laplacian $\Delta$ act only in the space variable.
\item If $f$ is a function and $\mu$ a measure, its convolution is $x\mapsto f\ast\mu(x):=\int f(x-y)\mu(\De y)$.
\eei


\subsection{The mean field Schr\"odinger problem and its equivalent formulations}\label{sec mfs intro}

We are given a so-called interaction potential $W:\RD\to\mathbb R,$  for which we assume 
\begin{align}\label{boundedhess}\tag{H1}
\text{$W$ is of class $\cC^{2}(\RD;\R)$ and symmetric, i.e. $W(\cdot)=W(-\cdot)$}, \\
\nonumber \sup_{z,v \in \mathbb{R}^d, |v|=1}  v\cdot \nabla^2W(z) \cdot v  < +\infty.
\end{align}
Besides the interaction potential, the data of the problem are a pair of probability measures $\mui,\muf$ on which we impose
\be\label{marginalhyp}\tag{H2}
\text{$\mu^{\mathrm{in}},\mu^{\mathrm{fin}}\in \cP_2(\RD)$ and $\absnent(\mu^{\mathrm{in}}),\absnent(\mu^{\mathrm{fin}})<+\infty$ },
\ee  
where the free energy or entropy functional $\absnent $ is defined for $\mu\in\cP_2(\RD)$ by

\be\label{eq def D}
\absnent(\mu) = \begin{cases} \int_{\RD} \log \mu(x) \mu(\De x)+ \int_{\RD} W\ast \mu(x) \mu(\De x), \quad & \mbox{if $\mu \ll \leb $}\\
+\infty, \quad & \mbox{otherwise.}
\end{cases}
\ee

In the above, and in the rest of the article, we shall make no distinction between a measure and its density against Lebesgue measure $\leb$, provided it exists. 

We recall that the McKean Vlasov dynamics is the non linear SDE 
\be\label{mckeanvlasovdynintro}
\begin{cases}
\De Y_t = - \nabla W\ast \mu_t (Y_t)\De t + \De B_t, \\
Y_0\sim \mu^{\mathrm{in}}, \quad \mu_t = \mathrm{Law}(Y_t), \quad \forall t\in[0,T].
\end{cases}
\ee
Under the hypothesis \eqref{boundedhess}, it is a classical result (see e.g. \cite[Thm 2.6]{cattiaux2008probabilistic})
that \eqref{mckeanvlasovdynintro} admits a unique strong solution whose law we denote $\rmP^{\text{\tiny{MKV}}}$. The functional $\absnent$ plays a crucial role in the sequel. For the moment, let us just remark that the marginal flow of the McKean-Vlasov dynamics may be viewed as the gradient flow of  $\frac{1}{2}\absnent$ in the Wasserstein space $(\cP_2(\RD),\mathcal W_2(\cdot,\cdot))$.

If $\rmP\in \cP_{1}(\Omega) $ is given, then the stochastic differential equation
$$
\left\{ 
\begin{array}{rrl}
     \De Z_t&=& -\nabla W \ast\rmP_t (Z_t)\De t+\De B_t , \\
     Z_0&\sim& \mu^{\mathrm{in}},
\end{array}
\right .
$$
admits a unique strong solution (cf.\ Section \ref{sec McKVl formulation proofs}) whose law we denote $\Gamma(\rmP)$. With this we can now introduce the main object of study of the article:
\begin{mydef}\label{def MFSP}
The mean field Schr\"odinger problem is 
\be\label{nonlinearSP}\tag{MFSP}
\inf \left\{\cH(\rmP | \Gamma(\rmP) ) \,:\,\,  \rmP \in \cP_{1}(\Omega), \,\rmP_0=\mu^{\mathrm{in}},\, \rmP_T=\mu^{\mathrm{fin}}\right\}.
\ee
Its optimal value, denoted $\nonlinentcost(\mui,\muf)$, is called mean field entropic transportation cost. Its optimizers are called mean field Schr\"odinger bridges (MFSB).
\end{mydef}

\begin{remark}
Of course, the choice of the pair potential $W$ as interaction mechanism is a particular one. Thus \eqref{nonlinearSP} is by far not the only mean field Schr\"odinger problem of interest. Also, it would have been easy to include in the dynamics a confinement (single-site) potential. However, since one of the goals of this article is to understand the role of the pair potential, we preferred not to do that. Indeed, in many situations, the single site potential may be the one that determines the long time behavior of mean field Schr\"odinger bridges, thus obscuring the role of the pair potential.
\end{remark}

It is not difficult to provide existence of optimizers for \eqref{nonlinearSP}. In the classical case, uniqueness is an easy consequence of the convexity of the entropy functional. However the rate function $\cH(\rmP|\Gamma(\rmP))$ is not convex in general. Therefore, it is not clear whether or not \eqref{nonlinearSP} admits a unique solution.

\begin{prop}\label{existence}
Grant \eqref{boundedhess},\eqref{marginalhyp}. Then \eqref{nonlinearSP} admits at least an optimal solution.
\end{prop}

\subsubsection{Large deviations principle (LDP)}

We start by deriving the LDP interpretation of \eqref{nonlinearSP}. We consider for $N\in \mathbb N$ the interacting particle system
\be\label{partsys}
\left\{ 
\begin{array}{rl}
     \De X_t^{i,N}=&-\frac{1}{N} \sum_{k=1}^N\nabla W(X_t^{i,N}-X_t^{k,N})\De t+\De B^i_t  \\
     X_0^{i,N}\sim& \mu^{\mathrm{in}}\,\, i=1,\ldots,N
\end{array}
\right .
\ee
where $\{B^i:i=1,\ldots, N\}$ are independent Brownian motions and $\{X_0^{i,N}:i=1,\ldots, N\}$ are independent to each other and to the Brownian motions. The theory of stochastic differential equations guarantees the strong existence and uniqueness for the above interacting particle system under \eqref{boundedhess},\eqref{marginalhyp}. 
In the next theorem we obtain a LDP for the sequence of empirical path measures;
 in view of the classical results of \cite{DaGae87}, it is not surprising that the LDP holds. However, even the most recent works on large deviations for weakly interacting particle systems such as \cite{BuDuFi12}  do not seem to cover the setting and scope of Theorem \ref{teo ldp intro}. { Essentially, this is because in those references the LDPs are obtained for a topology that is weaker than the $\cW_1$-topology, that is what we need later on.}

\begin{theorem}\label{teo ldp intro} In addition to \eqref{boundedhess},\eqref{marginalhyp} assume that 
\begin{align}\label{Ainitial} \textstyle 
\int_{\mathbb R^d}\exp(r|x|)\mu^{\mathrm{in}}(\De x)<\infty\text{ for all }r>0.
\end{align} 
Then the sequence of empirical measures $$ \textstyle  \left\{\frac{1}{N}\sum_{i=1}^N\delta_{X^{i,N}} ; N\in\mathbb N \right\},$$ satisfies the LDP on $\cP_{1}(\Omega)$ equipped with the $\mathcal W_1$-topology, with good rate function given by
\be\label{good rate form}
\cP_{1}(\Omega)\ni \rmP\mapsto \scrI(\rmP):=\left\{
\begin{array}{ll}
     \cH (\rmP|\Gamma(\rmP))&,\rmP\ll \Gamma(\rmP),  \\
     +\infty& ,\text{otherwise}.
\end{array}
\right .
\ee
\end{theorem}

In fact we will prove in Section \ref{sec formulations proofs} a strengthened version of Theorem \ref{teo ldp intro} where the drift term is much more general. For this, we will follow Tanaka's elegant reasoning \cite{Ta84}. 

\begin{remark}
Having a rate function implies heuristically 
\bes  \textstyle  \mathrm{Prob} \Big[ \frac{1}{N} \sum_{i=1}^N \delta_{X^{i,N}_{\cdot}} \approx \rmP  \Big] \approx \exp(-N \scrI(\rmP)).   \ees
Hence Problem \eqref{nonlinearSP} has the desired interpretation of finding the most likely evolution of the particle system conditionally on the observations (when $N$ is very large). This is in fact a consequence of the Gibbs conditioning principle; see \cite{dembo2009large} and the references therein.
\end{remark}

\subsubsection{McKean-Vlasov control and Benamou-Brenier formulation}\label{McK V fomulation intro}

We now reinterpret the mean field Schr\"odinger problem \eqref{nonlinearSP} in terms of McKean-Vlasov stochastic control (also known as mean field control).

\begin{lemma}\label{correctionexistence intro}
Let $\rmP$ be admissible for \eqref{nonlinearSP}. There exists a predictable process $(\alpha^{\rmP}_t)_{t\in[0,T]}$ s.t.\
\be\label{finitent2} \textstyle 
\bbE_{\rmP} \left[ \int_{0}^T |\alpha^{\rmP}_t|^2 \De t \right] <+\infty 
\ee
and so that
\be\label{finitent1} \textstyle 
X_t - \int_{0}^t \left( -\nabla W \ast\rmP_s (X_s) + \alpha^{\rmP}_s\right) \,\De s
\ee
has law $R^{\mu^{\mathrm{in}}}$ under $\rmP$. {The problem \eqref{nonlinearSP} is equivalent to 
\be\label{equivSch} \textstyle 
\inf \left\{ \frac{1}{2}\bbE_{\rmP} \left[ \int_{0}^T |\alpha^{\rmP}_t|^2 \De t \right] \,:\,\, \rmP \in \cP_{1}(\Omega), \,\rmP_0=\mu^{\mathrm{in}},\, \rmP_T=\mu^{\mathrm{fin}},\, \alpha^{\rmP}\text{ as in \eqref{finitent1}} \right \},
\ee
as well as to 
\begin{align}\label{equivSchSch}
\begin{split}
\inf  & \textstyle  \,\,\frac{1}{2}\bbE_{\rmP} \left[ \int_{0}^T |\Phi_t+\nabla W \ast\rmP_t (X_t) |^2 \De t \right] \\
 \text{s.t.}&  \textstyle \,\, \rmP \in \cP_{1}(\Omega), \,\rmP_0=\mu^{\mathrm{in}},\, \rmP_T=\mu^{\mathrm{fin}},\, \rmP\circ \left (X_\cdot-\int_0^\cdot\Phi_s\De s \right )^{-1}=R^{\mu^{\mathrm{in}}}.
\end{split}
\end{align}
}
\end{lemma}
The formulations \eqref{equivSch}-\eqref{equivSchSch} can be seen as McKean-Vlasov stochastic control problems. In the first case one is steering through $\alpha^{\rmP}$ part of the drift of a McKean-Vlasov SDE. In the second case one is controlling the drift $\Phi$ of a standard SDE but the optimization cost depends non-linearly on the law of the controlled process. In both cases, the condition $\rmP_T=\mu^{\mathrm{fin}}$ is rather unconventional. By analogy with the theory of mean field games, one could refer to \eqref{equivSch}-\eqref{equivSchSch} as \emph{planning McKean-Vlasov stochastic control problems}, owing to this type of terminal condition. 


{The third and last formulation of \eqref{nonlinearSP} we propose relates to the well known fluid dynamics representation of the Monge Kantorovich distance due to Benamou and Brenier (cf.\ \cite{villani2008optimal}) that has been recently extended to the standard entropic transportation cost \cite{gentil2015analogy,chen2016relation}. The interest of this formula is twofold: on the one hand it clearly shows that \eqref{nonlinearSP} is equivalent and gives a rigorous meaning to some of the generalized Schr\"odinger problems formally introduced in \cite{gentil2018dynamical,leger2019hopf}. On the other hand, it allows to interpret \eqref{nonlinearSP} as a control problem in the Riemannian manifold of optimal transport. This viewpoint, that we shall explore in more detail in Section \ref{OTFBSDEsection}, provides with a strong guideline towards the study of the long time behavior of Schr\"odinger bridges. }

We define the set $\mathcal A$ as the collection of all absolutely continuous curves $(\mu_t)_{t\in[0,T]}\subset \cP_2(\RD)$ (cf.\ Section \ref{OTsection}) such that $\mu_0=\mui,\mu_T=\muf$ and 
\begin{align*}
(t,z)\mapsto\nabla \log \mu_t(z) &\in L^2(\De\mu_t\De t),\\
(t,z)\mapsto\nabla W \ast \mu_t(z) &\in L^2(\De\mu_t\De t).
\end{align*}
We then define
\be\label{eq cost BB intro} \textstyle 
\nonlinentcost^{BB}(\mui,\muf):=\inf_{\substack{(\mu_t)_{t\in[0,T]}\in\mathcal A, \\ \partial_t \mu_t + \nabla \cdot (w_t \mu_t )=0}} \, \frac{1}{2}\int_0^T\int_{\RD} \left |w_t(z)+  \frac{1}{2}\nabla \log \mu_t(z) + \nabla W \ast \mu_t(z)\right|^2\mu_t(\De z) \De t.
\ee

\begin{theorem}\label{Teo BB intro}
Let \eqref{boundedhess},\eqref{marginalhyp} hold. Then $$\nonlinentcost(\mui,\muf) = \nonlinentcost^{BB}(\mui,\muf).$$ If $\rmP$ is optimal for \eqref{nonlinearSP} and the latter is finite, then $(\rmP_t)_{t\in[0,T]}$ is optimal in \eqref{eq cost BB intro} and its associated tangent vector field $w$ is given by $$ \textstyle  - \nabla W\ast \rmP_t(z)+ \Psi_t(z) -\frac{1}{2} \nabla\log \rmP_t ,$$ 
where $\Psi$ is as in Theorem \ref{martingalecorrection intro} below.

Conversely, if  $(\mu_t)_{t\in[0,T]}$ is optimal for $\nonlinentcost^{BB}(\mui,\muf)$ and the latter is finite, then there exists an optimizer of $\nonlinentcost(\mui,\muf) $ whose marginal flow equals $(\mu_t)_{t\in[0,T]}$. 
\end{theorem}

\subsection{Mean field Schr\"odinger bridges}
Leveraging the stochastic control interpretation, and building on the stochastic calculus of variations perspective, we obtain the following necessary optimality conditions for \eqref{nonlinearSP}.

\begin{theorem}\label{martingalecorrection intro}
Assume \eqref{boundedhess},\eqref{marginalhyp} and let $\rmP$ be optimal for \eqref{nonlinearSP}. Then there exist $\Psi \in \rmH_{-1}((\rmP_t)_{t\in[0,T]})$ such that 
\be\label{Markovopt1}
(\De t\times \De\rmP \text{-a.s.}) \quad \alpha^{\rmP}_t= \Psi_t(X_t),
\ee
where $(\alpha^{\rmP}_t)_{t\in [0,T]}$ is related to $\rmP$ as in Lemma \ref{correctionexistence intro}. The process $t\mapsto \Psi_t(X_t)$ is continuous\footnote{More precisely, it has a continuous version adapted to the $\rmP$-augmented canonical filtration.} and the process $(M_t)_{t\in[0,T]}$ defined by
\be\label{martingalecorrection1intro} \textstyle 
M_t:=\Psi_t(X_t) - \int_{0}^t \tilde{\bbE}_{\tilde{\rmP}} \left[ \nabla^2 W(X_s-\tilde{X}_s)\cdot(\Psi_s(X_s) - \Psi_s(\tilde{X}_s) )\right] \, \De s
\ee
is a continuous martingale under $\rmP$ on $[0,T[$, where $(\tilde{X}_t)_{t\in[0,T]}$ is an independent copy of $(X_t)_{t\in[0,T]}$ defined on some probability space $(\tilde{\Omega},\tilde{\cF},\tilde{\rmP})$ and $\tilde{\bbE}_{\tilde{\rmP}}$ denotes the expectation on $(\tilde{\Omega},\tilde{\cF},\tilde{\rmP)}$.
\end{theorem}

We shall refer to $\Psi$ as the \emph{corrector} of $\rmP$. Correctors will play an important role in the ergodic results. In this part, we give an interpretation of Theorem \ref{correctionexistence intro} in terms of stochastic analysis (FBSDEs) and partial differential equations.

\subsubsection{Planning McKean-Vlasov FBSDE for MFSB }\label{sec:FBSDE}

We consider the following McKean Vlasov forward-backward stochastic differential equation (FBSDE) in the unknowns $(X,Y,Z)$: 

\begin{equation}\label{FBSDE1intro}
\begin{cases}
\De X_t= -\tilde{\bbE}[\nabla W (X_t-\tilde{X}_t) ]\De t + Y_t\De t+ \De B_t\\
\De Y_t = \tilde{\bbE}\big[\nabla^2 W(X_t-\tilde{X}_t) \cdot(Y_t-\tilde{Y}_t)\big]\De t + Z_t\cdot \De B_t\\
X_0\sim\mu^{\mathrm{in}},\, X_T\sim\mu^{\mathrm{fin}}.
\end{cases}
\end{equation}
 As in the stochastic control interpretation of the mean field Schr\"odinger problem, here too the terminal condition $X_T\sim\mu^{\mathrm{fin}}$ is somewhat unconventional. We hence call this forward-backward system the \emph{planning McKean-Vlasov  FBSDE}.

Thanks to the results in Section \ref{McK V fomulation intro} we can actually solve \eqref{FBSDE1intro}. If $\rmP$ is optimal for \eqref{nonlinearSP} with associated $\Psi$ as recalled in Theorem \ref{martingalecorrection intro} above, all we need to do is take $Y_t:=\Psi_t(X_t)$ and reinterpret \eqref{finitent1} for the dynamics of the canonical process $X$ and \eqref{martingalecorrection1intro} for the dynamics of $Y$ (in the latter case using martingale representation). 

One remarkable aspect of this connection between Schr\"odinger problems and FBSDEs is that one can prove existence of solutions to such FBSDEs by a purely variational method. Indeed, we remark that \eqref{FBSDE1intro} is beyond the scope of existing FBSDE theory, such as Carmona and Delarue's \cite[Theorem 5.1]{CaDe15}. Further, we also obtained for free an extra bit of information: the constructed process $Y$ lives in $\rmH_{-1}((\rmP_t)_{t\in[0,T]})$. This is in tandem with the usual heuristic relating FBSDEs and PDEs (where $Y$ is conjectured to be an actual gradient) as explained in Carmona and Delarue's \cite[Remark 3.1]{CaDe13}. In fact, if we make the additional assumption that $Y_t=\nabla \psi_t(X_t)$  for some time-dependent sufficiently smooth potential $\psi_t(x)$, and we set $\mu_t=(X_{t})_{\#}\rmP$, then the first equation in \eqref{FBSDE1intro} tells that
\bes\label{FBSDE5intro} \textstyle 
\partial_t \mu_t - \frac{1}{2}\Delta \mu_t+ \nabla \cdot ((-\nabla W \ast \mu_t+\nabla \psi_t)\mu_t)=0,
\ees
whereas the second equation of \eqref{FBSDE1intro} 
%
and a standard argument with It\^o formula yields that for each $1\leq i \leq d$ that $\partial_{x_i} \psi_t(x)$ satisfies the PDE
\bes\label{FBSDE4intro} \textstyle 
(\partial_t + \mathcal{L}) \partial_{x_i} \psi_t(x) -  \int_{\RD} \nabla^2W(x-\tilde{x})\cdot(\partial_{x_i} \psi_s(x)-\partial_{x_i} \psi_s(\tilde{x})) \mu_t(\De \tilde{x}) =0,
\ees
where $\mathcal{L}$ is the non-linear generator $ \frac{1}{2}\Delta  - \big(\nabla  W\ast \mu_t-\nabla \psi_t \big) \cdot \nabla $. 

All in all we arrive at the PDE system\footnote{The Laplacian of a vectorial function is defined coordinate-wise.}:
\be\label{FBSDE6intro}
\begin{cases}
\partial_t \mu_t(x) - \frac{1}{2}\Delta \mu_t(x)+ \nabla \cdot ((-\nabla W \ast \mu_t(x)+\nabla \psi_t(x))\mu_t(x))=0\\
\partial_t \nabla \psi_t(x)+ \frac{1}{2}\nabla\Delta \psi_t(x)+\nabla^2 \psi_t(x) \cdot \big( -\nabla W \ast \mu_t(x) + \nabla \psi_t(x)\big)\\
= \int_{\RD} \nabla^2W(x-\tilde{x})\cdot(\nabla \psi_t(x)-\nabla \psi_t(\tilde{x})) \mu_t(\De \tilde{x}) ,\\
\mu_0(x) = \mui(x), \mu_T(x) = \muf(x).
\end{cases}
\ee
{It is worth stressing that the link between \eqref{FBSDE1intro} and \eqref{FBSDE6intro} can be established if the FBSDE solution $Y_t$ is a gradient vector field depending only on $t$ and $X_t$. We have gathered preliminary evidence that \eqref{FBSDE1intro} admits non Markovian solutions even in the simple case when $W=0$. This is in contrast with what is expected for standard FBSDEs \cite[Lemma 3.5]{CaDe13} whose boundary conditions are not of planning type.}

\subsubsection{Schr\"odinger potentials and the mean field planning PDE system}

The PDE system \eqref{FBSDE6intro} is the literal translation of the planning McKean-Vlasov FBSDE in the case when the process $Y$ is an actual gradient, $Y=\nabla\psi$. 
In the next corollary we show that if this is the case, and if $\psi$ is sufficiently regular, then \eqref{FBSDE6intro} can be rewritten as a system of two coupled PDEs, the first being a Hamilton-Jacobi-Bellman equation for $\psi$, and the second one being a Fokker-Planck equation. This type of PDE system is the prototype of a planning mean field game \cite{lasry2007mean}. 

\begin{cor}\label{mean field PDE sys}
Let $\rmP$ be an optimizer for \eqref{nonlinearSP}, $\Psi_{\cdot}(\cdot)$ be as in Theorem \ref{martingalecorrection intro} and set $\mu_t=\rmP_t$ for all $t\in[0,T]$. If $\mu_{\cdot}(\cdot)$ is everywhere positive and of class $\cC^{1,2}([0,T]\times \RD;\R)$ and $\Psi_{\cdot}(\cdot)$ is of class $\cC^{1,2}([0,T]\times \RD;\RD)$ then there exists $\psi:[0,T]\times\RD\rightarrow\R$ such that $\Psi_t(x) =\nabla \psi_t(x)$ for all $(t,x)\in[0,T]\times \RD$. Moreover, $(\psi_{\cdot}(\cdot),\mu_{\cdot}(\cdot))$ form a classical solution of
\begin{equation}\label{FBSDE master}
\begin{cases}
\partial_t \psi_t(x) + \frac{1}{2}\Delta \psi_t(x) + \frac{1}{2}|\nabla \psi_t(x)|^2 = \int_{\RD} \nabla W(x-\tilde{x}) \cdot( \nabla \psi_t(x)-\nabla \psi_t(\tilde{x}))\mu_t(\De \tilde{x}),\\
\partial_t \mu_t(x) - \frac{1}{2} \Delta \mu_t(x) + \nabla \cdot((-\nabla W \ast\mu_t(x) + \nabla \psi_t(x))\mu_t(x) )=0,\\
\mu_0(x) = \mui(x), \mu_T(x) = \muf(x)
\end{cases}
\end{equation}
\end{cor}

A fundamental result \cite{FOLL88} concerning the structure of optimizers in the classical Schr\"odinger problem is that their density takes a product form, i.e. 
\bes
\mu_t = \exp(\psi_t+\varphi_t),
\ees
where $\varphi_t(x)$, $\psi_t(x)$ solve respectively the forward and backward Hamilton Jacobi Bellman equation

\be\label{classical fg}
\begin{cases}
\,\,\,\,\partial_t \psi + \frac{1}{2}\Delta \psi + \frac{1}{2}|\nabla \psi|^2=0,\\
-\partial_t \varphi + \frac{1}{2}\Delta \varphi + \frac{1}{2}|\nabla \varphi|^2 =0.\\
\end{cases}
\ee

It is interesting to see that this structure is preserved in \eqref{nonlinearSP}, at least formally. The effect of having considered interacting Brownian particles instead of independent ones is reflected in the fact that the two Hamilton Jacobi Bellman PDEs are coupled not only through the boundary conditions but also through their dynamics.

\begin{cor}\label{cor Schr pot}
Using the same notation and under the same hypotheses of Corollary \ref{mean field PDE sys}, if we define  $\varphi:[0,T]\times \RD \rightarrow\R$ via 
\bes  \mu_t = \exp(-2 \nabla W\ast \mu_t + \varphi_t +\psi_t )\ees 
then $(\psi_{\cdot}(\cdot),\varphi_{\cdot}(\cdot))$ solves
\begin{equation*}\label{Schrodinger potentials}
    \begin{cases}
\,\,\,\,\partial_t \psi_t(x) + \frac{1}{2}\Delta \psi_t(x) + \frac{1}{2}|\nabla \psi_t(x)|^2 = \int \nabla W(x-\tilde{x}) \cdot( \nabla \psi_t(x)-\nabla \psi_t(\tilde{x}))\mu_t(\De\tilde{x}),\\
-\partial_t \varphi_t(x) + \frac{1}{2}\Delta \varphi_t(x) + \frac{1}{2}|\nabla \varphi_t(x)|^2 = \int \nabla W(x-\tilde{x}) \cdot( \nabla \varphi_t(x)-\nabla \varphi_t(\tilde{x}))\mu_t(\De\tilde{x}) .\\
\end{cases}
\end{equation*}
\end{cor}

\begin{remark}
In the article \cite{leger2019hopf} the authors implement another approach in order to generalize the decomposition \eqref{classical fg} to interacting systems based on the Hopf-Cole transform, whereas ours relies on time reversal. In particular, in their Example 4, they deal with a generalized Schr\"odinger problem that would correspond to erasing the term $\frac{1}{2}\nabla \log \mu_t$ in \eqref{eq cost BB intro}. It appears that the variables $\eta,\eta^*$ produced by the Hopf-Cole transformation do not coincide with the variables $\psi,\varphi$ we use here. Therefore these two approaches may be seen as complementary.
\end{remark}
%
%
%

\subsection{Convergence to equilibrium and functional inequalities}\label{sub quantitative intro}

Our aim is to show that MFSBs spend most of their time in a small neighborhood of the equilibrium configuration $\mu_{\infty}$, to study their long time behavior, and to derive a new class of functional inequalities involving the mean field entropic cost $\nonlinentcost(\mui,\muf)$.

\label{sec exp eq functional ineq}
 Throughout this section we make the assumption that $W$ is uniformly convex, ie.\ that
\be\label{W convexity ass}\tag{H3}
\exists \curv>0 \quad s.t. \quad  \forall z\in\RD, \quad \nabla^2 W(z) \geq \kappa \mathbb{I}_{d\times d},
\ee
where the inequality above has to be understood as an inequality between quadratic forms. Under \eqref{W convexity ass} the McKean Vlasov dynamics associated with the particle system \eqref{partsys} converges in the limit as $T\rightarrow +\infty$ to an equilibrium measure $\mu_{\infty}$, that is found by minimizing the functional $\absnent$ 
over the elements of $\cP_2(\RD)$ whose mean is the same as $\mui$. Existence and uniqueness of $\mu_{\infty}$ has been proven in \cite{mccann1997convexity}.

We shall often assume that $\mui$ and $\muf$ have the same mean:
 \be\label{same mean ass}\tag{H4} \textstyle 
\int_{\RD} x\, \mui (\De x) =  \int_{\RD} x\, \muf (\De x) .
 \ee
 
\begin{remark}
Assumption \eqref{W convexity ass} is a classical one ensuring exponential convergence rates for the McKean-Vlasov dynamics. It may be weakened in various ways, see the work \cite{carrillo2003kinetic} by Carrillo, McCann and Villani or the more recent \cite{bolley2013uniform} by Bolley, Gentil and Guillin, for instance. It is an interesting question to determine which of the results of this section still hold in the more general setup. Hypothesis \eqref{same mean ass} can be easily removed using the fact that the mean evolves linearly along any Schr\"odinger bridge (see Lemma \ref{linearbar} below); the reason for imposing \eqref{same mean ass} is that it allows to simplify the presentation of our results. We insist that all our quantitative results below use the technical assumptions \eqref{boundedhess} and \eqref{marginalhyp} only in order to make sure that all objects we are dealing with are well defined. If one can ensure this by other means, the results are still valid under the key assumption \eqref{W convexity ass}. 
 \end{remark}

\paragraph{Long time behavior of mean field games:} The articles \cite{cardaliaguet2012long,cardaliaguet2013long,cardaliaguet2013KAM,cardaliaguet2017long} study the asymptotic behaviour of dynamic mean field games showing convergence towards an ergodic mean field game with exponential rates. Following \cite{lasry2007mean}, we can associate to \eqref{FBSDE master} the ergodic problem with unknowns $(\lambda,\psi,\mu)$
\be \label{ergodic system}
\begin{cases}
\lambda + \frac{1}{2}\Delta \psi(x)+\frac{1}{2}|\nabla \psi(x)|^2 = \int_{\RD} \nabla W(x-\tilde x)\cdot(\nabla \psi(x)-\nabla \psi(\tilde{x}))\mu(\De \tilde{x}) \\
-\frac{1}{2}\Delta \psi(x) + \nabla \cdot ( (-\nabla W \ast \mu(x)+ \nabla \psi(x) )\mu(x)  )=0.
\end{cases}
\ee
This PDE system  expresses optimality conditions for the ergodic control problem corresponding to \eqref{equivSch}. It is easy to see that $(0,0,\mu_{\infty})$ is a solution of \eqref{ergodic system}. Therefore, we are addressing the same questions studied in the above mentioned articles. However, the equations we are looking at are quite different. A fundamental difference is that the coupling terms in \eqref{equivSch} are not monotone in the sense of \cite[Eq.(7) pg. 8]{cardaliaguet2013notes}.

\subsubsection{Exponential convergence to equilibrium and the turnpike property}

A key step towards the forthcoming quantitative estimates is to consider the time-reversed version of our mean field Schr\"odinger problem. For $\rmQ\in \cP(\Omega)$ the time reversal $\hat{\rmQ}$ is the law of the time reversed process $(X_{T-t})_{t\in[0,T]}$. In Lemma \ref{timerevoptimality} we prove that if $\rmP$ is an optimizer for \eqref{nonlinearSP}, then $\hat{\rmP}$ optimizes
\bea\label{timerevnonlinearSP}
\inf \left\{\cH(\rmQ | \Gamma(\rmQ) ) \,:\,  \rmQ \in \cP_{1}(\Omega),\,  \rmQ_0=\mu^{\mathrm{fin}},\, \rmQ_T=\mu^{\mathrm{in}} \right \}.
\eea

The optimality of $\hat{\rmP}$ implies the existence of an associated process $\hat\Psi$ as described in Theorem \ref{martingalecorrection intro}. We show at Theorem  \ref{consquantity1intro} below that the function
\be\label{consquantity2intro}
[0,T]\ni t\mapsto \bbE_{\rmP}[\Psi_t(X_t)\cdot\hat{\Psi}_{T-t}(\hat{X}_{T-t})]
\ee
 is a constant, that we denote $\cq(\mu^{\mathrm{in}},\mu^{\mathrm{fin}})$ and call the \emph{conserved quantity}. Naturally this quantity depends also on $T$ but we omit this from the notation. 


Theorem \ref{exp conv eq} confirms the intuition that mean field Schr\"odinger bridges are localized around $\mu_{\infty}$ providing an explicit upper bound for $\nent(\rmP_{t})$ along any MFSB, where
\be\label{def relative nent}
\nent(\mu) = \absnent(\mu)-\absnent(\mu_{\infty}).
\ee
We recall that $\mu_{\infty}$ is found by minimizing $\absnent$ among all elements of $\cP_2(\RD)$ whose mean is the same as $\mu$. If $\absnent$ is thought of as a free energy, then $\nent$ should be thought of as a divergence (from equilibrium).\footnote{In principle the equilibrium measure $\mu_{\infty}$ associated to $\rmP_t$ varies with $t$. However, because of Lemma \ref{linearbar} below and assumption \eqref{same mean ass}, we obtain that $\mu_{\infty}$ does not depend on $t$ (and neither on $T$) once $\mui,\muf$ are fixed.}

\begin{theorem}\label{exp conv eq}
Assume \eqref{boundedhess}-\eqref{same mean ass} and let $\rmP$ be an optimizer for \eqref{nonlinearSP}. For all $t\in[0,T]$ we have
\be \label{consquantconvexity1intro} \textstyle \nent(\rmP_t)\leq \frac{\sinh(2\curv (T-t))}{\sinh(2\curv T)}\Big(\nent(\mu^{\mathrm{in}})-\frac{\cq(\mu^{\mathrm{in}},\mu^{\mathrm{fin}})}{2\curv} \Big)
 +\frac{\sinh(2\curv t)}{\sinh(2\curv T)}\Big( \nent(\mu^{\mathrm{fin}})-\frac{\cq(\mu^{\mathrm{in}},\mu^{\mathrm{fin}})}{2\curv}\Big) + \frac{\cq(\mu^{\mathrm{in}},\mu^{\mathrm{fin}})}{2\curv}.\ee
 Moreover, for all fixed $\theta\in(0,1)$ there exists a decreasing function $B(\cdot)$  such that
 \be\label{exp entropy bound} \nent(\rmP_{\theta T})\leq B(\curv)(\nent(\mui) +\nent(\muf)) \exp(-2 \curv \min\{\theta, 1-\theta \} T) \ee
 uniformly in $T \geq 1$. 
\end{theorem}

In particular, since $\nent(\rmP_{\theta T})$ dominates $\mathcal W_2(\rmP_{\theta T},\mu_{\infty})$ (see e.g. \cite[(ii), Thm 2.2 1]{carrillo2003kinetic}), we obtain that $\rmP_{\theta T}$ converges exponentially fast to $\mu_{\infty}$ and that the exponential rate is proportional to $\curv$. {The proof of \eqref{consquantconvexity1intro} is done by bounding the second derivative of the function $t\mapsto \nent(\rmP_t)$ along Schr\"odinger bridges with the help of the logarithmic Sobolev inequality established in \cite{carrillo2003kinetic}. To obtain \eqref{exp entropy bound} from \eqref{consquantconvexity1intro} we use a  functional inequality for the conserved quantity and a Talagrand inequality for $\nonlinentcost(\mui,\muf)$, that are the content of Theorem \ref{consquantity1intro} and Corollary \ref{talagrand} below. It is worth mentioning that the estimates \eqref{consquantconvexity1intro},\eqref{exp entropy bound} (as well as \eqref{talagrandconsqtyintro} below) appear to be new even for the classical Schr\"odinger bridge problem and have not been anticipated by the heuristic articles \cite{gentil2018dynamical,leger2019hopf}. Conversely, the above mentioned estimates admit a geometrical interpretation in the framework of Otto calculus that allows to formally extend their validity to the whole class of problems studied in \cite{gentil2018dynamical}.}

\begin{remark}
The exponential rate in \eqref{exp entropy bound} has a sharp dependence on $\curv$. To see this, fix $\mui$ and choose $\muf=\mkv_T$. Then it is easy to see that the restriction of $\rmP^{\text{\tiny{MKV}}}$ to the interval $[0,T]$ is an optimizer for \eqref{nonlinearSP}. Setting $\theta=1/2$ and considering \eqref{exp entropy bound} for $T=2t$ we arrive at \bes \forall t\geq 1/2, \quad \nent(\mkv_t) \leq B(\curv) \exp(-2\curv t) \ees
Thus, we obtain the same exponential rate  as in \cite{carrillo2003kinetic}\footnote{{Some doubt on the numeric value of the exponential rates may arise from the fact that in our definition of $\tilde{\cF}$, there is no $1/2$ in front of $W$, as it is the case in \cite[Eq. 1.3]{carrillo2003kinetic}. However, as we pointed out before, the McKean-Vlasov dynamics for the particle system \eqref{partsys intro} is the gradient flow of $1/2\tilde{\cF}$ and not of $\tilde{\cF}$.}}, that is easily seen to be optimal under the assumption that $W$ is $\curv$-convex. A similar argument can be used to show the optimal dependence of the rate in $\theta$.
\end{remark}
In the previous theorem we showed that, when looking at a timescale that is proportional to $T$, the marginal distribution of any Schr\"odinger bridge is exponentially close to $\mu_{\infty}$. Here we show that for a \emph{fixed} value of $t$, we have an exponential convergence towards the law of the McKean-Vlasov dynamics $\rmP^{\text{\tiny{MKV}}}$, see \eqref{mckeanvlasovdynintro}.

\begin{theorem}\label{thm exp mkv intro}
Assume \eqref{boundedhess}-\eqref{same mean ass} and let $\rmP$ be an optimizer for \eqref{nonlinearSP}. For all $t\in[0,T]$ we have
\be\label{correctionbound2}  
\mathcal W_2^2(\rmP_t,\rmP^{\text{\tiny{MKV}}}_t) \leq 2t \left(\, \frac{\nent(\mu^{\mathrm{in}})}{\exp(2\curv T)-1} + \frac{\exp(2\curv T)-\exp(2\curv (T-t))}{  \exp(2\curv (T-t))-1} \frac{\nent(\mu^{\mathrm{fin}})}{\exp(2\curv T)-1} \right)
\ee
\end{theorem}

In particular, the above theorem tells that $\mathcal W_2^2(\rmP_t,\rmP^{\text{\tiny{MKV}}}_t)$ decays asymptotically at least as fast as $\exp(-2\curv T)$ when $T$ is large. 

\subsubsection{Functional inequalities for the mean field entropic cost}

It is well known that analysing the evolution of entropy-like functionals along the so-called displacement interpolation of optimal transport has far reaching consequences in terms functional inequalities \cite{von2005transport}. Since \eqref{nonlinearSP} provides with an alternative way of interpolating between probability measures, it is tempting to see if it leads to new functional inequalities involving the cost $\nonlinentcost(\mui,\muf)$. Here, we present a Talagrand and an HWI inequality that we used in order to study the long time behavior of MFSBs. They generalize their respective counterparts in \cite{talagrand1996transportation},\cite{otto2000generalization}. Both inequalities are based on another upper bound for the evolution of $\nent$ along MFSBs that we present now. This bound generalizes to the mean field setup the results of \cite{conforti2018second}: 

\begin{theorem}\label{nonlinearentropybound1intro}
Assume \eqref{boundedhess}-\eqref{same mean ass} and let $T>0$ be fixed. If $\rmP$ is an optimizer for \eqref{nonlinearSP}, then for all $t\in [0,T]$ we have
\bea\label{nonlinearentropyboundintro}
\nonumber \nent(\rmP_t) \leq \frac{\exp(2\curv (T-t))-1}{\exp(2\curv T)-1}\nent(\mu^{\mathrm{in}}) +  \frac{\exp(2\curv T)-\exp(2\curv (T-t))}{\exp(2\curv T)-1} \nent(\mu^{\mathrm{fin}})\\ - \frac{ (\exp(2\curv (T-t) )-1) (\exp(2\curv t )-1 )}{\exp(2\curv T)-1} \nonlinentcost(\mu^{\mathrm{in}},\mu^{\mathrm{fin}}).
\eea
\end{theorem}

The above result may be seen as a rigorous version of some of the heuristic arguments put forward in \cite{gentil2018dynamical} and \cite{leger2019hopf}, upon slightly modifying the definition of $\nonlinentcost$.

The following Talagrand inequality, obtained as a corollary of Theorem \ref{nonlinearentropybound1intro}, tells that the mean field entropic cost grows at most linearly with $\nent$.

\begin{cor}[A Talagrand Inequality]\label{talagrand}
Assume \eqref{boundedhess}-\eqref{same mean ass}.  Then for all $T>0$ we have 
\be\label{talagrand1intro}
\forall t\in(0,T),\quad \nonlinentcost(\mu^{\mathrm{in}},\mu^{\mathrm{fin}}) \leq \frac{1}{\exp(2\curv t)-1} \nent(\mu^{\mathrm{in}}) + \frac{\exp(2\curv(T-t))}{\exp(2\curv(T-t))-1}\nent(\mu^{\mathrm{fin}}).
\ee
In particular, choosing $\muf=\mu_{\infty}$ leads to
\be\label{talagrand2intro}
 \nonlinentcost(\mu^{\mathrm{in}},\mu_{\infty}) \leq \frac{1}{\exp(2\curv T)-1} \nent(\mu^{\mathrm{in}}).
\ee
\end{cor}

Theorem \ref{nonlinearentropybound1intro} also yields an entropic HWI inequality where the Wasserstein distance is replaced by  the conserved quantity $\cq$ in the first term on the rhs and by the mean field entropic cost in the second term. An extra positive contribution $\frac{1}{4}\nonlinfish$ is present in the first term. Our interpretation is that this compensates for the fact that in the ``gain" term we put the cost $\nonlinentcost$, that is larger than the squared Wasserstein distance. 
In order to state the HWI inequality, we introduce the non linear Fisher information functional $\nonlinfish$ defined for $\mu\in\cP_2(\RD)$ by

\be\label{non linear fisher}
\nonlinfish(\mu) = \begin{cases} \int_{\RD} \Big|\nabla \log \mu + 2\nabla W \ast \mu (x)\Big|^2 \mu(\De x), \quad & \mbox{if $\nabla \log \mu \in L^2_{\mu}$ }\\
+\infty \quad & \mbox{otherwise.}
\end{cases}
\ee
where by $\nabla \log \mu \in L^2_{\mu}$ we mean $\mu\ll\leb$ and that $\log \mu$ is an absolutely continuous function on $\RD$ whose derivative is in $L^2_{\mu}$. The non linear Fisher information can be seen to be equal to the derivative of the free energy $\absnent$ along the marginal flow of the McKean Vlasov dynamics.

\begin{cor}[An HWI Inequality]\label{HWI} Assume \eqref{boundedhess}-\eqref{W convexity ass} and choose $\muf=\mu_{\infty}$. If $\rmP$ is an  optimizer for \eqref{nonlinearSP}  and $t\mapsto \nonlinfish(\rmP_t)$ is continuous in a right neighbourhood of $0$, then
\be\label{HWI consquant}
\nent(\mui) \leq \frac{1-\exp(-2\curv T)}{2\curv} \left(\nonlinfish(\mui) \left( \frac{1}{4}\nonlinfish(\mui) - \cq(\mui,\mu_{\infty}) \right)\right)^{1/2} - (1-\exp(-2\curv T)) \nonlinentcost(\mui,\mu_{\infty}). 
\ee
\end{cor}

It is worth noticing that by letting $T\rightarrow +\infty$ in the above HWI inequality we obtain the logarithmic Sobolev inequality \cite[Thm 2.2]{carrillo2003kinetic}. Indeed, $\nonlinentcost(\mui,\mu_{\infty})$ is always non negative and we shall see at Thoerem \ref{consquantity1intro} below that $\cq(\mui,\mu_{\infty}) \rightarrow 0$. The short time regime is also interesting. Indeed, if $W=0$, $\nonlinentcost(\mui,\mu_{\infty})$ is the standard entropic cost and we have under suitable hypothesis on $\mui$ (see \cite{Mik04})
\be\label{short time convergence cost} \textstyle 
\lim_{T\rightarrow0} T \nonlinentcost(\mui,\mu_{\infty})= \frac{1}{2} \mathcal W_2^2(\mui,\mu_{\infty}).
\ee
The heuristic arguments put forward in \cite{gentil2018dynamical}  tell that \eqref{short time convergence cost} is expected to be true even when $W$ is a general potential satisfying \eqref{boundedhess}. Following again \eqref{short time convergence cost}, one also expects that 
\be\label{short time convergence consqty}  \textstyle \lim_{T \rightarrow 0}  \frac{T^2}{4}\nonlinfish(\mui) - T^2 \cq(\mui,\mu_{\infty}) =  \mathcal  W^2_2(\mui,\mu_{\infty}).
\ee
Putting \eqref{short time convergence cost} and \eqref{short time convergence consqty} together we obtain an heuristic justification of the fact that in the limit as $T\rightarrow 0$ \eqref{HWI consquant} becomes the classical HWI inequality
put forward in \cite{carrillo2003kinetic}, namely
\bes
\nent(\mui) \leq \mathcal W_2(\mui,\mu_{\infty}) \nonlinfish(\mui)^{1/2} - \curv \mathcal W^2_2(\mui,\mu_{\infty}).
\ees

\begin{remark}
In Corollary \ref{HWI} we asked $t\mapsto \nonlinfish(\rmP_t)$ to be continuous in a right neighbourhood of $0$. This is somewhat annoying. However, we were not able conclude that \eqref{boundedhess} and \eqref{marginalhyp} imply the desired continuity, although we could establish that they grant the continuity of $\nonlinfish(\rmP_t)$ on any open subinterval of $[0,T]$.
\end{remark}
Our last result is a functional inequality that establishes a hierarchical relation between the conserved quantity and the mean field entropic cost: the former is exponentially small in $T$ and $\curv$ in comparison with the latter. {We may refer to this as an \emph{energy-transport} inequality since the conserved quantity may be geometrically interpreted as the total energy of a physical system (cf \cite[Corollary 1.1]{conforti2018second}).}
\begin{theorem}\label{consquantity1intro}
Assume \eqref{boundedhess}-\eqref{same mean ass} and let $\rmP$ be an optimizer. Then the function
\be\label{consquantity2intro}
[0,T]\ni t\mapsto \bbE_{\rmP}[\Psi_t(X_t)\cdot\hat{\Psi}_{T-t}(\hat{X}_{T-t})]
\ee
is constant. Denoting this constant by $\cq(\mu^{\mathrm{in}},\mu^{\mathrm{fin}})$, we have
\be\label{talagrandconsqtyintro}
|\cq(\mu^{\mathrm{in}},\mu^{\mathrm{fin}})| \leq \frac{4\curv}{\exp(\curv T)-1} \left(\nonlinentcost(\mu^{\mathrm{in}},\mu^{\mathrm{fin}})\nonlinentcost(\mu^{\mathrm{fin}},\mu^{\mathrm{in}}) \right)^{1/2}.
\ee
\end{theorem}

\section{Connections with optimal transport}\label{OTFBSDEsection}

In this section we shall see how the results of this article relate to the Riemannian calculus on $\cP_2(\RD)$ introduced by Otto \cite{otto2001geometry}, at least formally. The link is rooted in a seemingly novel connection between (McKean-Vlasov) FBSDEs and second order ODEs in the Riemannian manifold of optimal transport that we find of independent interest. To better understand this connection, let us begin by recalling that in the seminal article \cite{jordan1998variational} it is proven that the marginal flow of the pathwise SDE
\be\label{ItoSDE} \textstyle 
\De X_t = - \nabla U(X_t)\De t +\De B_t
\ee
can be interpreted as the \emph{gradient flow} of the entropy functional 
$$ \textstyle  \mu \mapsto \frac{1}{2}\int_{\RD} \log \mu(x) \mu(\De x) +\int_{\RD} U(x)\mu(\De x)$$
w.r.t.\ the $2$-Wasserstein metric. Thus, first order It\^o SDEs provide with \emph{probabilistic representations} for first order ODEs in the Riemannian manifold of optimal transport. Of course, since a path measure is not fully determined by its one time marginals, the SDE \eqref{ItoSDE} contains more information than the gradient flow equation. It has been shown in \cite{conforti2018second} that the marginal flow of a classical Schr\"odinger bridge satisfies a second order ODE, more precisely a \emph{Newton's law} in which the acceleration field is the Wasserstein gradient of the Fisher information functional. The natural question is then: What \emph{pathwise} (second order) SDE governs the dynamics of a Schr\"odinger bridge and yields a probabilistic representation for the associated Newton's law? In order to answer this, let us first recall some notions of Otto calculus.


\subsection{Second order calculus on $\cP_2(\RD)$.} In the next lines, we sketch the ideas behind the Riemannian calculus on $\cP_2(\RD)$. It would be impossible to provide a self-contained introduction in this work and we refer to \cite{villani2008optimal} or \cite{gigli2012second} for detailed accounts. The main idea is to equip $\cP_2(\RD)$ with a Riemannian metric such that the associated geodesic distance is $\mathcal W_2(\cdot,\cdot)$. To do this, one begins by identifying the tangent space $\tansp$ at $\mu\in\cP_2(\RD)$ as the space closure in $L^2_{\mu}$ of the subspace of gradient vector fields 
\bes  \textstyle 
\tansp= \overline{\{ \nabla \varphi, \varphi \in \cC^{\infty}_c(\RD)\}}^{L^2_{\mu}}.
\ees
The velocity  (first derivative) of a sufficiently regular curve $[0,T]\ni t\mapsto \mu_t\in \cP_2(\RD)$ is then defined by looking at the only  solution $v_t(x)$ of the continuity equation
\bes  \textstyle 
\partial_t \mu_t + \nabla \cdot (v_t \mu_t)=0
\ees
such that $v_t\in\tanspt$ for all $t\in[0,T]$. Finally, the Riemannian metric (Otto metric) $\langle\cdot,\cdot \rangle_{\tansp}$ is defined by
\be\label{ottometric}
\langle \nabla \varphi,\nabla \psi \rangle_{\tansp} = \int_{\RD}\nabla \varphi\cdot \nabla \psi(x) \, \mu(\De x) .
\ee

It can be seen that the constant speed geodesic curves associated to the Riemannian metric we have introduced coincide with the displacement interpolations of optimal transport and that the corresponding geodesic distance is indeed $\mathcal W_2(\cdot,\cdot)$. This makes it possible to carry out several explicit calculations. In particular, we can compute the gradient $\otgrad \cF$ and the Hessian $\othess \nent $ of a smooth functional $\cF:\cP_2(\RD)\rightarrow\R$. At least formally, we have
\begin{align*}
\langle\otgrad \cF,\nabla \varphi\rangle_{\tansp}&=  \textstyle \frac{\De}{\De h} \cF((id+h \nabla \varphi)_{\#}\mu) \Big|_{h=0}\\
\langle \nabla \varphi, \othess_{\mu} \cF(\nabla \varphi)\rangle_{\tansp}& =  \textstyle  \frac{\De^2}{\De h^2} \cF((id+h \nabla \varphi)_{\#}\mu) \Big|_{h=0},
\end{align*}
where we used the notation $\#$ for the push forward. In particular, setting $W=0$ for simplicity in \eqref{non linear fisher} we obtain that the classical Fisher information functional $\cI$ has a gradient  that can be computed with the rules above. One obtains that (cf.\ \cite{von2012optimal})
\bes
\otgrad \cI(\mu) = -2\nabla \Delta  \log \mu - \nabla |\nabla \log \mu|^2.
\ees

The Levi-Civita connection associated to the Riemannian metric \eqref{ottometric} can also be explicitly computed with the help of the orthogonal projection operator $\Pi_{\mu}:L^2_{\mu} \rightarrow \tansp$. To do this, consider a regular curve $(\mu_t)_{t\in[0,T]}$ with velocity $(v_t)_{t\in[0,T]}$ and a tangent vector field $t\mapsto u_t\in \tanspt$ along $(\mu_t)_{t\in[0,T]}$. It turns out that if one defines the covariant derivative $\frac{\mathbf{D}}{\De t}u_t$ of $(u_t)_{t\in[0,T]}$ along $(\mu_t)_{t\in[0,T]}$ as the vector field
\bes  \textstyle 
\frac{\mathbf{D}}{\De t}u_t = \Pi_{\mu_t}\left(\partial_t u_t + 
\mathrm{D} u_t \cdot v_t \right)
\ees
 then this covariant derivative satisfies the compatibility with the metric and the torsion-free identity, i.e. it is the Levi-Civita connection. The acceleration of the curve $(\mu_t)_{t\in[0,T]}$ is then the covariant derivative of the velocity along the curve, i.e.

\be\label{newtonlaw}  \textstyle 
\frac{\mathbf{D}}{\De t}v_t = \partial_t v_t + \frac{1}{2}\nabla |v_t|^2.
\ee

\subsection{Newton's laws and FBSDEs}

 According to the above discussion the Newton's law in $(\cP_2(\RD), \langle.,.\rangle_{\cT_{\cdot}\cP_2})$
\be\label{OTnewton}
\begin{cases}\frac{\mathbf{D}}{\De t} v_t = \frac{1}{8} \otgrad  \cI (\mu_t)\\
\mu_{0}=\mui, \mu_T=\muf \end{cases}
\ee
provides with a geometrical interpretation for the PDE system (see \cite{conforti2018second} for more details)
\be\label{PDEnewton}
\begin{cases}
 \partial_t \mu_t(x) + \nabla \cdot (\nabla \phi_t(x) \mu_t(x))=0 \\
\partial_t \nabla \phi_t(x) + \frac{1}{2}\nabla |\nabla \phi_t(x)|^2 = -\frac{1}{4}\nabla \Delta \log \mu_t(x) - \frac{1}{8}\nabla |\log \mu_t(x)|^2 \\
\mu_0=\mu^{\mathrm{in}},\mu_T=\mu^{\mathrm{fin}},
\end{cases}
\ee
where to derive the latter equation we observe that the requirement that $v_t\in\tanspt$ for all $t\in[0,T]$ is formally equivalent to $v_t=\nabla \phi_t $ for some time dependent potential $(t,x)\mapsto \phi_t(x)$.


As we have seen in section \ref{sec:FBSDE}, solutions of the FBSDE
\begin{equation}\label{FBSDE1otsection}
\begin{cases}
\De X_t=   Y_t\De t+ \De B_t\\
\De Y_t =  Z_t\cdot \De B_t\\
X_0\sim\mu^{\mathrm{in}},\, X_T\sim\mu^{\mathrm{fin}},
\end{cases}
\end{equation}
 having the additional property that  $Y_t=\nabla \psi_t(X_t)$ yield a probabilistic representation for
 
\be\label{FBSDE6otsection}
\begin{cases}
\partial_t \mu_t(x) - \frac{1}{2}\Delta \mu_t(x)+ \nabla \cdot (\nabla \psi_t(x)\mu_t(x))=0,\\
\partial_t \nabla \psi_t(x)+ \frac{1}{2}\nabla\Delta \psi_t(x)+\nabla^2 \psi_t(x) \cdot \nabla \psi_t(x)=0,\\
\mu_0(x) = \mui(x), \mu_T(x) = \muf(x).
\end{cases}
\ee
 Some tedious though standard calculations allow to see that the change of variable $\phi_t= -\frac{1}{2}\log \mu_t + \psi_t$ transforms the PDE system \eqref{FBSDE6otsection} in \eqref{PDEnewton}. Summing up, we have obtained the following
 \paragraph{Informal statement}
 \emph{We have:
 \bei 
 \item[(i)] If $(X_t,Y_t,Z_t)_{t\in[0,T]}$ is a solution for the FBSDE \eqref{FBSDE1otsection} such that $Y_t=\nabla\psi_t(X_t)$ for some time-varying potential $\psi$, then the marginal flow $(\mu_t)_{t\in[0,T]}$ of $X_t$ is a solution for the Newton's law \eqref{OTnewton}.
 \item[(ii)] If $\rmP$ is the (classical) Schr\"odinger bridge between $\mui$ and $\muf$, then under $\rmP$ the canonical process $(X_t)_{t\in[0,T]}$ is such that there exist processes $(Y_t)_{t\in[0,T]}$,$(Z_t)_{t\in[0,T]}$ with the property that $(X_t,Y_t,Z_t)_{t\in[0,T]}$ is a solution for \eqref{FBSDE1otsection} and $Y_t$ is as in (i)
\eei 
}
We leave it to future work to prove a rigorous version of the informal statement above. On the formal level, there is no conceptual difficulty in extending it  to include the interaction potential $W$. Essentially, the only difference is that one has to deal with the non linear Fisher information functional $\cI_{\nent}$ instead of $\cI$.

Beside its intrinsic interest, the parallelism between Newton's laws and FBSDEs is very useful when  studying the long time behavior of the latter. Indeed, the Riemannian structure underlying \eqref{OTnewton} allows to find tractable expressions for the first and second derivative of entropy-like functionals along the marginal flow of the FBSDE.
\begin{remark}
Classical Schr\"odinger bridges are $h-$transforms in the sense of Doob \cite{Doob1957}. Therefore, one can also describe their dynamics with a first order SDE and a PDE that encodes the evolution of the drift field. This is not strictly speaking a probabilistic representation of \eqref{OTnewton} since there is already a PDE involved. Our FBSDE approach may be viewed as a way to interpret in a pathwise sense the PDE governing the dirft in the $h-$transform representation.
\end{remark}

\section{The mean field Schr\"odinger problem and its equivalent formulations: proofs}\label{sec formulations proofs}

In this part we complement the discussion undertaken in Section \ref{sec mfs intro} and provide the proofs of the results stated therein. { This section is organized into four subsections so that

\bei 
\item Subsection \ref{sec LDP} contains the proof of a more general version of Theorem \ref{teo ldp}, that is made with the help of several lemmas, 
\item Subsection \ref{sec McKVl formulation proofs} is where we prove Proposition \ref{existence}, Lemma \ref{correctionexistence intro} and Theorem \ref{martingalecorrection intro}.
\item Theorem \ref{Teo BB intro} is proven in subsection \ref{sec ben bren proofs}.
\item Finally, Corollary \ref{mean field PDE sys} and \ref{cor Schr pot} are proven in Subsection \ref{sec schr pot}.
\eei
In the whole section, apart from subsection \ref{sec LDP} that has its own assumptions, we always assume that \eqref{boundedhess},\eqref{marginalhyp} are in force, even if we do not write them down explicitly in the statements of the lemmas and propositions.
}

\subsection{A large deviations principle for particles interacting through their drifts}
\label{sec LDP}

We consider for $N\in \mathbb N$ the interacting particle system
$$
\left\{ 
\begin{array}{rl}
     \De X_t^{i,N}=&\frac{1}{N} \sum_{k=1 }^Nb\left(t,X^{i,N},X^{k,N}\right)\De t+\De B^i_t  \\
     X_0^{i,N}\sim& \mui,\,\,\, i=1,\ldots,N.
\end{array}
\right .
$$
where $\{B^i:i=1,\ldots, N\}$ are independent Brownian motions and $\{X_0^{i,N}:i=1,\ldots, N\}$ are independent to each other and to the Brownian motions. Regarding the drift $b$, we assume
\begin{align}
&[0,T]\times \Omega\times\Omega \ni (t,\omega,\bar\omega)\mapsto b(t,\omega,\bar\omega)\in\mathbb R^d \text{ is progressively measurable,}\label{Adrift1}\\ 
&\textstyle |b(t,\omega^1,\bar\omega^1)-b(t,\omega^2,\bar\omega^2)|\leq C\left\{\sup_{s\leq t}|\omega^1_s-\omega_s^2|+ \sup_{s\leq t}|\bar\omega^1_s-\bar\omega_s^2|\right\}\label{Adrift2}\\
&\textstyle \int_0^T|b(s,0,0)|\De s\leq C,\label{Adrift3}
\end{align}
for some constant $C>0$ and all $(t,\omega^1,\omega^2,\bar\omega^1,\bar\omega^2)\in [0,T]\times \Omega^4$. Finally, regarding the measure $\mui$ we assume that 
\begin{align}\label{Ainitial} \textstyle
\int_{\mathbb R^d}\exp(r|x|^\beta)\mui(\De x)<\infty\text{ for all }r>0.
\end{align}
We stress that the usual theory of stochastic differential equations guarantees the strong existence and uniqueness for the above interacting particle system. Furthermore, if $\rmP\in \cP_{1}(\Omega) $ then the same arguments show that the stochastic differential equation
$$
\left\{ 
\begin{array}{rl}
     \De X_t^{\rmP}=& \left[\int b\left(t,X^{\rmP},\bar\omega\right) \rmP(\De \bar\omega)\right ]\De t+\De B_t  \\
     X_0^{\rmP}\sim& \mui,
\end{array}
\right .
$$
admits a unique strong solution. We denote $\Gamma(\rmP)$ the law of $X^\rmP$. We can now state the main result of this part, which contains Theorem \ref{teo ldp intro} as a very particular case.

\begin{theorem}\label{teo ldp}Let $\beta\in[1,2)$ and assume \eqref{Adrift1},\eqref{Adrift2},\eqref{Adrift3},\eqref{Ainitial}. Then the sequence of empirical measures $$ \textstyle \left\{\frac{1}{N}\sum_{i=1}^N\delta_{X^{i,N}_{\cdot}} :N\in\mathbb N \right\},$$ satisfies a LDP on $\cP_{\beta}(\Omega)$ equipped with the $\mathcal W_\beta$-topology, with good rate function given by
\be\label{good rate form b}
\cP_{\beta}(\Omega)\ni \rmP\mapsto \scrI(\rmP):=\left\{
\begin{array}{ll}
     \cH (\rmP|\Gamma(\rmP))&,\rmP\ll \Gamma(\rmP),  \\
     +\infty& ,\text{otherwise}.
\end{array}
\right .
\ee
\end{theorem}

The result is sharp, in that it fails for $\beta=2$; see \cite{WWW10}. We follow Tanaka's reasoning \cite{Ta84} in order to establish this large deviations result. We remark that the assumption on exponential moments \eqref{Ainitial} is only used in the proof of Theorem \ref{teo ldp}, and not in the results preceding this proof. 

For $\rmQ\in\cP_\beta(\Omega)$ we consider the equation
\be \textstyle 
Y_t(\omega)=\omega_t+\int_0^t\left[\int b(s,Y(\omega),Y(\bar\omega)) \rmQ(\De \bar\omega) \right]\De s.\label{eq path integral equation}
\ee

\begin{lemma}\label{lem well posed iter}
Take $Y^{(0)}_t(\omega):=\omega_0$, $\rmQ\in\cP_\beta(\Omega)$, and consider the iterations 
$$ \textstyle Y_t^{(n+1)}(\omega)=\omega_t+\int_0^t\left[\int b(s,Y^{(n)}(\omega),Y^{(n)}(\bar\omega)) \rmQ(\De \bar\omega) \right]\De s,\,\, s\leq T.$$
Then
\begin{itemize}
    \item[a)] The iteration is well-defined $\omega$-by-$\omega$ (in particular, the $\rmQ$-integrals are well-defined and finite) and in fact $\sup_n\bbE_\rmQ\left[ \sup_{t\leq T}|Y^{n}_t|^\beta \right]$ is finite.
    \item[b)]  For each $\omega\in\Omega$ the sequence $\{Y^{(n)}(\omega)\}_{n\in\N}$ is convergent in the sup-norm to some limiting continuous path $Y^{(\infty)}(\omega)$. Further $\bbE_\rmQ\left[ \sup_{t\leq T}|Y^{(\infty)}_t|^\beta \right]<\infty$, $\bbE_\rmQ\left[ \sup_{t\leq T}|Y^{(\infty)}_t-Y^{(n)}_t| \right]\to 0$, and $Y^{(\infty)}$ is adapted to the canonical filtration.
\end{itemize}
\end{lemma}

\begin{proof}
From the Lipschitz assumption on $b$ we first derive
\begin{align} \textstyle 
    \sup_{s\leq t}|Y^{(n+1)}_s|\leq \sup_{s\leq t}|\omega_s|+\int_0^T|b(s,0,0)|\De s +C\int_0^t\sup_{r\leq s}|Y^{(n)}_r|\De r +C\int_0^t\bbE_\rmQ\left[\sup_{r\leq s}|Y^{(n)}_r|\right ]\De r.
\end{align}
Raising this to $\beta$, taking expectations and using Jensen's inequality, we derive
$$  \textstyle \bbE_\rmQ\left[\sup_{s\leq t}|Y^{(n+1)}_s|^\beta\right ]\leq C'\left (1+ \bbE_\rmQ\left[\sup_{s\leq T}|\omega_s|^\beta\right ] +\int_0^t\bbE_\rmQ\left[\sup_{r\leq s}|Y^{(n)}_r|^\beta\right ]\De r \right ),$$
where $C'$ only depends on $T$ and $\beta$. From this we establish for some $R\geq 0$ that $$ \textstyle \sup_n \bbE_\rmQ\left[\sup_{s\leq t}|Y^{(n)}_s|^\beta\right ] \leq Re^{Rt}.$$
Now denote $\Delta^{n}_t:=\sup_{s\leq t}|Y^{(n)}_s-Y^{(n-1)}_s|$. Again by the Lipschitz property
$$ \textstyle \Delta^{n+1}_t\leq C\int_0^t \left\{\Delta^{n}_s + \bbE_\rmQ [\Delta^{n}_s]  \right\} \De s,$$
which we can bootstrap to obtain 
$$ \textstyle \Delta^{n+1}_t+\bbE_\rmQ [\Delta^{n+1}_t] \leq 3C\int_0^t \left\{\Delta^{n}_s + \bbE_\rmQ [\Delta^{n}_s]  \right\} \De s.$$ {Observe that $\Delta_t^1\leq 2\sup_{s\leq T}|\omega_s-\omega_0|+C$, so from the above inequality we obtain by induction that $\Delta^{n+1}_t+\bbE_\rmQ [\Delta^{n+1}_t] \leq C'' \frac{t^n}{n!}$.}
From this $\{\Delta^{n}_T+\bbE_\rmQ [\Delta^{n}_T]\}_{n\in\N}$ is (for each $\omega$) summable in $n$, so the same happens to $\{\Delta^{n}_T\}_{n\in\N}$ and therefore the uniform limit of the $Y^{(n)}$ exists for all $\omega$. We denote by $Y^{(\infty)}$ this limit. By Fatou's lemma $\bbE_\rmQ\left[ \sup_{t\leq T}|Y^{(\infty)}_t|^\beta \right]<\infty$. Since $\left(\bbE_\rmQ [\Delta^{n}_T]\right)_{n\in\N}$ is summable we must also have $\bbE_\rmQ\left[ \sup_{t\leq T}|Y^{(\infty)}_t-Y^{(n)}_t| \right]\to 0$. Since clearly each $Y^{(n)}$ is adapted so is $Y^{(\infty)}$ too.
\end{proof}

\begin{lemma}\label{lem unique pointwise}
For any $\rmQ\in\cP_{\beta}$ there exists a unique adapted continuous process satisfying \eqref{eq path integral equation} pointwise. Denoting $Y^\rmQ$ this process, we further have
$$\rmQ\circ (Y^\rmQ)^{-1}\in\cP_\beta(\Omega).$$
\end{lemma}

\begin{proof}
If $X$ and $Y$ are solutions, then the Lipschitz assumption on $b$ implies
$$ \textstyle \bbE_\rmQ\left[\sup_{s\leq t}|Y_s-X_s| \right]\leq K \int_0^t \bbE_\rmQ\left[\sup_{r\leq s}|Y_r-X_r| \right]\De r,$$
so from Gr\"onwall we derive $\bbE_\rmQ\left[\sup_{s\leq T}|Y_s-X_s| \right]=0$. With this, and using again the Lipschitz assumption on $b$, we find $$ \textstyle \sup_{s\leq t}|Y_s-X_s| \leq K \int_0^t \sup_{r\leq s}|Y_r-X_r| \De r,$$
so by Gr\"onwall we deduce that $X=Y$ pointwise. For the existence of a solution we employ Point (b) of Lemma \ref{lem well posed iter}, taking limits in the iterations therein (the exchange of limit and integral is justified by the Lipschitz property of $b$). Finally $\rmQ\circ (Y^\rmQ)^{-1}\in\cP_\beta(\Omega)$ follows by Point (b) of Lemma \ref{lem well posed iter} too.
\end{proof}

Thanks to this result we can define the operator
\begin{align}\label{defi op Theta}
\begin{split}
    \Theta:& (\cP_\beta,\mathcal W_\beta) \to (\cP_\beta,\mathcal W_\beta)\\
    &\,\,\,\,\rmQ\mapsto\Theta(\rmQ):=\rmQ\circ (Y^\rmQ)^{-1},
    \end{split}
\end{align}
where $Y^\rmQ$ denotes the unique solution of \eqref{eq path integral equation}.

\begin{lemma}\label{lem identification}
$Y^{\rmR^{\mui}}$ is the unique strong solution to the McKean-Vlasov SDE
$$
\left\{ 
\begin{array}{rl}
     \De Z_t=& \left[\int b\left(t,Z,\bar\omega\right) \rmP(\De \bar\omega)\right ]\De t+\De B_t  \\
     Z\sim& \rmP,\, Z_0\sim \mui.
\end{array}
\right .
$$
Furthermore, if $\{X^{i,N}:i\leq N, N\in\mathbb N\}$ is the aforementioned interacting particle system, which is driven by 
$\{B^i:i\in\mathbb N\}$ independent Brownian motions stared like $\mui$, then
\be \label{eq push forward empirical} \textstyle 
\Theta\Big(\frac{1}{N}\sum_{i=1}^N\delta_{B^i_{\cdot}}  \Big )= \frac{1}{N}\sum_{i=1}^N\delta_{X^{i,N}_{\cdot}},\,\, a.s.  
\ee
\end{lemma}

\begin{proof}
That $Y^{\rmR^{\mui}}$ is a solution to the McKean-Vlasov SDE is clear since $\omega$ is a Brownian motion under $\rmR^{\mui}$. That the solution is unique follows by observing that the drift in this SDE is Lipschitz jointly in $Z$ and $P=\text{Law}(Z)$, from where usual arguments apply. For the second point, consider first $\omega^1,\dots,\omega^N$ continuous paths and define $ \rmQ=\frac{1}{N}\sum_{i=1}^N\delta_{\omega^i}$. Then for all $1\leq i\leq N$ we have
$$ \textstyle Y^\rmQ_t(\omega^i)=\omega^i_t + \int_0^t\Big (\frac{1}{N}\sum_{k\leq N}b(s,Y^\rmQ(\omega^i),Y^\rmQ(\omega^k))  \Big )\De s.$$
Replacing the deterministic paths $\omega^1,\dots,\omega^N$ by those of $B^1,\dots,B^N$ we conclude.
\end{proof}

The key observation is that $\frac{1}{N}\sum_{i=1}^N\delta_{B^i} $ satisfies a large deviations principle in $\cP_\beta(\Omega)$ equipped with the $\mathcal W_\beta$ topology, with good rate function given by the relative entropy $\cH (\cdot|\rmR^{\mui})$. This is true for $\beta<2$ under our exponential moments assumption \eqref{Ainitial}, but fails for $\beta=2$, as follows easily from \cite{WWW10}. By Lemma \ref{lem identification} we may derive, via the contraction principle (\cite[Theorem 4.2.1]{dembo2009large}) a large deviations principle for 
$$ \textstyle \Big\{\frac{1}{N}\sum_{i=1}^N\delta_{X^{i,N}} :N\in\mathbb N \Big\},$$
if we could only establish the continuity of $\Theta$. This is our next step.

\begin{lemma}\label{lem op Lip}
$\Theta$ is Lipschitz-continuous and injective.
\end{lemma}

\begin{proof}
We first prove the Lipschitz property. Let $\pi$ be a coupling with first marginal $\rmQ$ and second marginal $\rmP$. Denoting $(\omega,\bar\omega)$ the canonical process on $\Omega\times\Omega$, and by the Lipschitz assumption on $b$, we have
$$ \textstyle \bbE_\pi\left [ \sup_{s\leq t}|Y_s^\rmQ(\omega)-Y_s^\rmP(\bar\omega)|^\beta\right ]\leq K\int_0^t \bbE_\pi\left [ |Y_s^\rmQ(\omega)-Y_s^\rmP(\bar\omega)|^\beta \right]\De s  +\bbE_\pi\left [ \sup_{s\leq t}|\omega_s-\bar\omega_s|^\beta\right ].$$
By Gr\"onwall we have 
$$ \textstyle \bbE_\pi\left [ \sup_{s\leq T}|Y_s^\rmQ(\omega)-Y_s^\rmP(\bar\omega)|^\beta\right ]\leq K'  \bbE_\pi\left [ \sup_{s\leq T}|\omega_s-\bar\omega_s|^\beta\right ],$$
so taking infimum over such $\pi$ we conclude that $$\mathcal W_\beta(\Theta(\rmQ),\Theta(\rmP))\leq K' \mathcal W_\beta(\rmQ,\rmP) .$$

We now prove that $\Theta$ is injective. Let $\rmP=\Theta(\rmQ)=\Theta(\hat\rmQ)$. By definition we have $\rmQ$-a.s.
\begin{align*}
    \omega_t&= \textstyle Y_t^\rmQ(\omega)-\int_0^t\left[ \int b(s,Y_s^\rmQ(\omega),Y_s^\rmQ(\bar\omega)) \rmQ(\De\bar\omega)\right]\De s = Y_t^\rmQ(\omega)-\int_0^t\left[ \int b(s,Y_s^\rmQ(\omega),\bar\omega) \rmP(\De\bar\omega)\right]\De s,
\end{align*}
and the same holds for $\hat\rmQ$ instead of $\rmQ$. Denoting $$ \textstyle F(\omega):=\omega-\int_0^\cdot\left[ \int b(s,\omega,\bar\omega) \rmP(\De\bar\omega)\right]\De s,$$
we therefore have
\begin{align*}
    \omega_t&= F(Y^\rmQ)_t\,\, (\rmQ-a.s.),\\
    \omega_t&= F(Y^{\hat\rmQ})_t\,\, (\hat\rmQ-a.s.).
\end{align*}
Hence $\rmQ=\Theta(\rmQ)\circ(F)^{-1}=\rmP\circ(F)^{-1}=\Theta(\hat\rmQ)\circ(F)^{-1}=\hat\rmQ$.
\end{proof}

We can now provide the proof of Theorem \ref{teo ldp}:

\begin{proof}[Proof of Theorem \ref{teo ldp}]
As we have observed, if $\{B^i:i\in\mathbb N\}$ is and iid sequence of $\rmR^{\mui}$-distributed processes, then $\frac{1}{N}\sum_{i=1 }^N\delta_{B^i_{\cdot}} $ satisfies a large deviations principle in $\cP_\beta(\Omega)$ equipped with the $\mathcal W_\beta$ topology, with good rate function given by the relative entropy $\cH (\cdot|\rmR^{\mui})$. By \eqref{eq push forward empirical}, and since $\Theta:(\cP_\beta,\mathcal W_\beta) \to(\cP_\beta,\mathcal W_\beta)$ is continuous, the contraction principle establishes that 
$\{\frac{1}{N}\sum_{i=1}^N\delta_{X^{i,N}_{\cdot}} :N\in\mathbb N \}$  satisfies a large deviations principle in $\cP_\beta(\Omega)$ equipped with the $\mathcal W_\beta$ topology. Since $\Theta$ is injective the good rate function is given by
$$
\tilde \scrI(\rmP):=\left\{
\begin{array}{ll}
     \cH(\Theta^{-1}(\rmP)|\rmR^{\mui})& \text{ if $\rmP\in$ range$(\Theta)$} \\
     +\infty &\text{ otherwise}.
\end{array}
\right .
$$
{In fact observe that if $\rmP\in$ range$(\Theta)$ and $\Theta^{-1}(\rmP)\ll\rmR^{\mui}$ then\footnote{Let $\rmP=\Theta(\rmQ)$ for $\rmQ\ll\rmR^{\mui}$. The process $Y^\rmQ$ satisfies pointwise $\De Y_t^\rmQ=\De \omega_t+\bar{b}(t,Y^\rmQ)\De t $, where $\bar{b}(t,y)=\int b(t,y,Y^\rmQ(\bar{\omega}))\rmQ(\De\bar{\omega})$. We have $\rmR^{\mui}\circ(Y^\rmQ)^{-1}\ll \rmR^{\mui}$ since in fact their relative entropy is finite. Hence, if $\rmR^{\mui}(A)=0$ then $\rmR^{\mui}((Y^\rmQ)^{-1}(A))=0$, and so $\rmQ\ll\rmR^{\mui}$ implies $\rmQ((Y^\rmQ)^{-1}(A))=0$ therefore $\rmP(A)=0$ as desired. } $\rmP\ll\rmR^{\mui}$, so
$$
\tilde \scrI(\rmP):=\left\{
\begin{array}{ll}
     \cH(\Theta^{-1}(\rmP)|\rmR^{\mui})& \text{ if $\rmP\in$ range$(\Theta)$ and $\rmP\ll\rmR^{\mui}$} \\
     +\infty &\text{ otherwise}.
\end{array}
\right .
$$
Now take $\rmP\in$ range$(\Theta)$ and call $\rmQ=\Theta^{-1}(\rmP)$. It is immediate by the definition of $\Gamma(\cdot)$ that $\Gamma(\rmP)=\rmR^{\mui}\circ(Y^\rmQ)^{-1}$. On the other hand observe that the filtration generated by $Y^\rmQ$ is equal to the canonical filtration: indeed $Y^\rmQ$ is adapted and conversely
$$ \textstyle \omega_t=Y^\rmQ_t-\int_0^t \left[\int b(s,Y^\rmQ_s,\bar\omega)\rmP(\De\bar\omega) \right ] \De s=:h_t(Y^Q),$$
so the canonical process is $Y^\rmQ$-adapted. From this $$ \textstyle \frac{\De\Big(\rmQ\circ(Y^\rmQ)^{-1}\Big)}{\De\Big(\rmR^{\mui}\circ(Y^\rmQ)^{-1}\Big)}=\bbE_{\rmR^{\mui}}\left[\frac{\De\rmQ}{\De\rmR^{\mui}}|\sigma(Y^\rmQ)\right]=\frac{\De\rmQ}{\De\rmR^{\mui}}\circ h.$$ Hence
$$ \textstyle \cH(\rmP|\Gamma(\rmP))=\cH(\rmQ\circ(Y^\rmQ)^{-1}|\rmR^{\mui}\circ(Y^\rmQ)^{-1})=\bbE_{\rmQ\circ(Y^\rmQ)^{-1}}\left[ \log \frac{\De\rmQ}{\De\rmR^{\mui}}\circ h  \right ]= \cH(\rmQ|\rmR^{\mui})=\cH(\Theta^{-1}(\rmP)|\rmR^{\mui}),$$
and therefore
$$\tilde \scrI(\rmP)=\left\{
\begin{array}{ll}
     \cH(\rmP|\Gamma(\rmP))& \text{ if $\rmP\in$ range$(\Theta)$ and $\rmP\ll\rmR^{\mui}$} \\
     +\infty &\text{ otherwise}.
\end{array}
\right .
$$
The next step is to show that $\rmP\ll\rmR^{\mui}$ 
implies $\rmP\in$ range$(\Theta)$. In fact, denote by $\tau$ the adapted transformation
$$ \textstyle \omega\mapsto \tau_t(\omega)=\omega_t-\int_0^t\int b(s,\omega,\bar\omega)\rmP(\De \bar\omega)\De s.$$
On the other hand call $X^{\rmP}$ the unique adapted pointwise solution to
$$ \textstyle X_t^{\rmP}=\omega_0+\int_0^t\left[\int b\left(s,X^{\rmP},\bar\omega\right) \rmP(\De \bar\omega)\right ]\De s + \omega_t,$$
which exists by Lemma \ref{lem unique pointwise} applied to the drift $\int b(\cdot,\cdot,\bar{\omega})\rmP(\De \bar{\omega})$.
As we recall in Lemma \ref{shiftinvert}, $X^{\rmP}$ and $\tau$ are $\rmP$-a.s.\ inverses if $\rmP\ll\rmR^{\mui}$, since the above drift is Lipschitz. Now introduce $\rmQ:=\rmP\circ (\tau)^{-1} $, so that $\rmQ\circ(X^{\rmP})^{-1}=\rmP$ and in particular 
$$ \textstyle X_t^{\rmP}=\omega_0+\int_0^t\left[\int b\left(s,X^{\rmP},X^{\rmP}(\bar\omega)\right) \rmQ(\De \bar\omega)\right ]\De s + \omega_t.$$
By Lemma \ref{lem unique pointwise} we have $\Theta(\rmQ):=Q\circ(Y^{\rmQ})^{-1}=Q\circ(X^{\rmP})^{-1} =\rmP$.
}

{
We have arrived at 
$$\tilde \scrI(\rmP)=\left\{
\begin{array}{ll}
     \cH(\rmP|\Gamma(\rmP))& \text{ if $\rmP\ll\rmR^{\mui}$} \\
     +\infty &\text{ otherwise}.
\end{array}
\right .
$$
To obtain the desired form \eqref{good rate form b} of the rate function it suffices to use Lemma \ref{entropyequivalence} in the Appendix.
}
\end{proof}

\subsection{McKean-Vlasov formulation and planning McKean-Vlasov FBSDE}\label{sec McKVl formulation proofs}

%

\paragraph{Proof of Lemma \ref{correctionexistence intro} and Proposition \ref{existence}}

Under \eqref{boundedhess} for any $\rmP\in \cP_{1}(\Omega)$ the vector field $$[0,T]\times\RD \ni(t,x)\mapsto -\nabla W\ast \rmP_t(x):=-\int_{\RD}\nabla W(x-z)\rmP_t(\De z), $$ is very well-behaved. Precisely: 

\begin{lemma}\label{driftregularity}
Let $\rmP \in \cP_{1}(\Omega)$ and grant \eqref{boundedhess}. Then the time-dependent vector field $(t,x)\mapsto -\nabla W\ast \rmP_t(x) $ belongs $ \cC^{0,1}([0,T]\times\RD;\RD)$ and is uniformly Lipschitz in the space variable.
\end{lemma}

\begin{proof}
We begin by proving continuity. Fix $t,x$ and $(t_n,x_n) \rightarrow (t,x)$. The sequence  $\nabla W(x_n-X_{t_n})$ converges pointwise to $\nabla W(x-X_{t})$, since $X$ is the (continuous) canonical process. By the fundamental theorem of calculus and \eqref{boundedhess} we have $|\nabla W(x_n-X_{t_n})|\leq C_1+C_2\sup_{s\in[0,T]}|X_s|$.  Since $\rmP\in \cP_{1}(\Omega)$, we may use dominated convergence to conclude $\bbE_{\rmP}[\nabla W(x_n-X_{t_n})] \rightarrow \bbE_{\rmP}\left[\nabla W(x-X_{t})\right]$. The space Lipschitzianity of $-\nabla W \ast \rmP_t$ follows from \eqref{boundedhess}. Space differentiability follows similarly from \eqref{boundedhess} and dominated convergence. 
\end{proof}

We will often make use of the next technical lemma, whose proof we defer to the appendix:


\begin{lemma}\label{finitent}
Let $\mu\in \cP_2(\RD)$ and $\bar b $ be of class $\cC^{0,1}([0,T] \times \RD;\RD)$ and such that
\be\label{lipschitz}
\forall t\in[0,T], \, x,y\in\RD \quad |\bar b(t,x)-\bar b(t,y)|\leq C|x-y|
\ee

for some $C<+\infty$. Define $\bar{\rmR}$ as the law of the SDE 
\be\label{SDELip}
\De X_t = \bar b(t,X_t) \De t + \De B_t, \quad X_0 \sim \mu
\ee

and let $\rmP \in \cP(\Omega)$ with $X_0 \sim \mu$. The following are equivalent

\bei 
\item[(i)] $\cH(\rmP|\bar{\rmR})<+\infty$.
\item[(ii)] There exist a $\rmP$-a.s.\ defined adapted process $(\bar{\alpha}_t)_{t\in[0,T]}$ such that
\be\label{finitent9} \textstyle 
\bbE_{\rmP} \left[ \int_{0}^T |\bar{\alpha}_t|^2 \De t \right] <+\infty 
\ee
and 
\be\label{finitent12} \textstyle 
X_t - \int_{0}^t[ \bar b(s,X_s) + \bar{\alpha}_s]\,\De s
\ee
is a Brownian motion under $\rmP$.
\eei 

Moreover, under condition (ii)   
\be\label{finitent13} \textstyle  \cH(\rmP | \bar{\rmR} ) = \frac{1}{2} \bbE_{\rmP}\left[ \int_{0}^T |\bar{\alpha}_t|^2 \De t \right] \ee
and 
\be\label{secondmoment} \textstyle 
\bbE_{\rmP}\left[\,\sup_{t\in[0,T]} |X_t|^2+|\bar b(t,X_t)^2|\,\right]<+\infty.
\ee
In particular under either (i) or (ii) we have $\rmP\in\cP_2(\Omega)$.
\end{lemma}
We turn to proving Lemma \ref{correctionexistence intro} stated in the introduction:
\begin{proof}[Proof of Lemma \ref{correctionexistence intro}]
Define the vector field $\bar b(t,z):= - \nabla W \ast \rmP_t (z)$. Lemma \ref{driftregularity} grants that $\bar b$ fulfills the hypotheses of Lemma \ref{finitent}, giving the desired conclusions.
\end{proof}

We can prove Proposition \ref{existence} of the introduction, concerning the existence of MFSBs. Recall the definition of $\Gamma(\rmP)$ and \eqref{nonlinearSP} from the introduction.

\begin{proof}[Proof of Proposition \ref{existence}]
Let $\rmR^{\mui}$ be the law of the Brownian motion started at $\mui$. \eqref{marginalhyp} grants that the classical Schr\"odinger problem (namely wrt.\ Brownian motion) is admissible. To see this, it suffices to verify that the coupling $\mui \otimes \muf$ is admissible for the static version of the Schr\"odinger problem \cite[Def 2.2]{LeoSch} and then use the equivalence between the static and dynamic versions \cite[Prop 2.3]{LeoSch}. Therefore, there exist some $\rmP \in \cP(\Omega)$ such that $\rmP_0=\mui$ and $\cH(\rmP|\rmR^{\mui})<+\infty$. Lemma \ref{finitent} (or its specialization Lemma \ref{wienerent} in the appendix) yields that $\rmP\in \cP_{1}(\Omega)$. On the other hand Lemma \ref{entropyequivalence} in the appendix proves that for any $\rmP\in\cP(\Omega)$ $\cH(\rmP|\Gamma(\rmP))<+\infty$ if and only if $\cH(\rmP|\rmR^{\mui})<+\infty$. Thus \eqref{nonlinearSP} is admissible as well. Now observe that $\rmP\mapsto \cH(\rmP|\Gamma(\rmP))$ is lower semicontinuous in $\cP_\beta(\Omega)$, since on the one hand the relative entropy is jointly lower semicontinuous in the weak topology, and on the other hand $\Gamma$ is readily seen to be continuous in $\cP_{1}(\Omega)$. Recalling the definition of the operator $\Theta$ given in \eqref{defi op Theta}, to finish the proof we only need to justify that $$\theta_M:=\{\rmP\in\cP_1(\Omega):\,\cH(\Theta^{-1}(\rmP)|\rmR^{\mui})\leq M,\, \rmP_0=\mui \},$$
is relatively compact in $\cP_1(\Omega)$ for each $M$, since the proof of Theorem \ref{teo ldp} established\footnote{This part of the proof did not use the existence of exponential moments for $\mui$. If we assume existence of exponential moments, then the compactness of $\theta_M$ follows from Theorem \ref{teo ldp}, since the rate function must be good.} that $\cH(\Theta^{-1}(\rmP)|\rmR^{\mui})= \cH(\rmP|\Gamma(\rmP))$ if $\rmP\ll\rmR^{\mui}$. Now remark that
$$\theta_M\subset \Theta\left( \{\rmQ:\cH(\rmQ|\rmR^{\mui})\leq M, \rmQ_0=\mui\}\right)\subset \Theta\left( \{\rmQ:\cH(\rmQ|\rmR^{\gamma})\leq \bar M, \rmQ_0=\mui\}\right),$$ 
since by the decomposition of the entropy we have
\bes  \textstyle 
\cH(\rmP|\rmR^{\gamma}) = \cH(\mui|\gamma)+ \cH(\rmP|\rmR^{\mui}),
\ees
and by Assumption \eqref{marginalhyp} $$ \textstyle \cH(\mui|\gamma)=\int\log\mui(x)\mui(\De x)-\int \log(\gamma(x))\mui(\De x)=\int\log\mui(x)\mui(\De x)+c-\int \frac{|x|^2}{2}\mui(\De x) <\infty.$$
As $\Theta$ is per Lemma \ref{lem op Lip} Lipschitz in $\cP_1(\Omega)$, it remains to prove that $\{\cH(\rmQ|\rmR^{\gamma})\leq\bar M\}$ is $\cW_1$-compact. This can be easily done by hand, or by invoking Sanov Theorem in the $\cW_1$-topology for independent particles distributed according to $\rmR^\gamma$ (see e.g.\ \cite{WWW10}), finishing the proof.
\end{proof}




\paragraph{Proof of Theorem \ref{martingalecorrection intro}} We split the proof into two propositions, namely Propositions \ref{Markovopt} and \ref{martingalecorrection}.
We begin by addressing the issue of Markovianity of the minimizers. Recall the definition of $\rmH_{-1}((\mu_t)_{t\in[0,T]})$ given under `frequently used notation.' We rely strongly on the work \cite{cattiaux1995large} by Cattiaux and L\'eonard for the proof of the following result:
\begin{prop}\label{Markovopt}
Let $\rmP$ be optimal for \eqref{nonlinearSP}. Then there exists $\Psi \in \rmH_{-1}((\rmP_t)_{t\in[0,T]})$ such that 
\be\label{Markovopt1}
(\De t\times \De\rmP \text{-a.s.}) \quad \alpha^{\rmP}_t= \Psi_t(X_t),
\ee
where $(\alpha^{\rmP}_t)_{t\in [0,T]}$ is given in Lemma \ref{correctionexistence intro}.
\end{prop}

\begin{proof}
If $\rmP$ be optimal for \eqref{nonlinearSP}, then it is also optimal for
\be\label{constrainedSP}
\inf \left\{\cH(\rmQ | \Gamma(\rmP) ) : \rmQ \in \cP_{1}(\Omega),\, \rmQ_t = \rmP_t\text{ for all }t\in[0,T]\right\},
\ee
since $\Gamma(\rmP)$ only depends on the marginals of $\rmP$. The above problem is an instance of \cite{cattiaux1995large}, ie.\ its optimizer is a so-called critical Nelson process. However, the drift of the path-measure $\Gamma(\rmP)$ may not fulfill the hypotheses in \cite{cattiaux1995large}. For this reason we need to make a slight detour. Let $\theta^n\in \cC^{\infty}_{c}([0,T]\times \RD)$ and $\rmR^n$ be defined as in Lemma \ref{entropyapprox} in the appendix, meaning that $\nabla \theta^n_{\cdot}(\cdot)$ converges to $-\nabla W \ast \rmP_t(z)$ in $\rmH_{-1}((\rmP_t)_{t\in[0,T]})$ and that $\rmR^n$ is the law of
\bes
\De Y_t = \nabla \theta^n_t(Y_t)\De t + \De B_t, \quad Y_0 \sim \mui\in\cP_2(\mathbb R^d).
\ees
For any $n$ consider the problem 
\be\label{constrainedSPn}
\min \left\{\cH(\rmQ | \rmR^n ) :\rmQ \in \cP_{1}(\Omega), \, \rmQ_t = \rmP_t\text{ for all }t\in[0,T]\right\}.
\ee
Using \cite[Lemma 3.1,Theorem 3.6]{cattiaux1995large} we obtain that for all $n$ the unique optimizer $\bar{\rmP}$ of \eqref{constrainedSPn} is the same for all $n$, and is such that there exists $\Phi\in \rmH_{-1} ((\rmP_t)_{t\in[0,T]})$ such that 
\be\label{markovcorrection}  \textstyle 
X_t - \int_{0}^t \Phi_s(X_s) \De s
\ee
is a Brownian motion under $\bar{\rmP}$. Lemma \ref{driftregularity} grants that if we set  $\bar b(t,z)= -\nabla W\ast\rmP_t(z)$ then the hypotheses of Lemma \ref{finitent} are met.  Since $\cH(\rmP|\Gamma(\rmP))<+\infty$, we derive from \eqref{secondmoment} therein that
\bes  \textstyle 
\bbE_{\bar{\rmP}} \left[ \int_{0}^T | \nabla W\ast \rmP_t(X_t)|^2\De t\right]=\bbE_{\rmP} \left[ \int_{0}^T | \nabla W\ast \rmP_t(X_t)|^2\De t\right]<+\infty.
\ees
Hence
\be\label{finitent14} \textstyle 
\bbE_{\bar{\rmP}} \left[ \int_{0}^T |\Phi_t(X_t)+ \nabla W\ast \rmP_t(X_t)|^2\De t\right]<+\infty.
\ee
Using the implication $(ii)\Rightarrow(i)$ of Lemma \ref{finitent} we finally obtain that $\cH(\bar{\rmP}|\Gamma(\rmP))<+\infty$ and therefore that we can use Lemma \ref{entropyapprox} for the choice $\rmQ=\bar{\rmP}$ therein. 

Now consider $\rmQ$ admissible for \eqref{constrainedSP} and such that $\cH(\rmQ|\Gamma(\rmP))<+\infty$. Using Lemma \ref{entropyapprox} twice we obtain
\bes
\cH(\bar{\rmP}|\Gamma(\rmP))= \liminf_{n\rightarrow +\infty} \cH(\bar{\rmP}|\rmR^n) \leq \liminf_{n\rightarrow +\infty} \cH(\rmQ|\rmR^n)=\cH(\rmQ|\Gamma(\rmP))
\ees
Thus $\bar{\rmP}$ is also an optimizer for \eqref{constrainedSP}. But then $\bar{\rmP}=\rmP$ since \eqref{constrainedSP} can have at most one minimizer by strict convexity of the entropy and convexity of the admissible region. Combining \eqref{markovcorrection} with \eqref{finitent1} we get that $\int_{0}^t \left( -\nabla W \ast\rmP_s (X_s) + \alpha^{\rmP}_s - \Phi_s(X_s)\right) \De s$
is a continuous martingale with finite variation. But then it is constant $\rmP$-a.s. The conclusion follows setting $\Psi_t(z):= \Phi_t(z)+\nabla W\ast \rmP_t(z)$ and observing that $\nabla W\ast \rmP_{\cdot}(\cdot)\in\rmH_{-1}((\rmP_t)_{r\in[0,T]})$.
\end{proof}

Notice that the above proposition proves the first half of Theorem \ref{martingalecorrection intro} from the introduction. We now establish the second half of this result:

\begin{prop}\label{martingalecorrection}
Assume that $\rmP$ is optimal for \eqref{nonlinearSP}. {Then $\Psi_t(X_t)$ has a continuous version adapted to the $\rmP$-augmented canonical filtration,} and the process $(M_t)_{t\in[0,T]}$ defined by
\be\label{martingalecorrection1} \textstyle 
M_t:=\Psi_t(X_t) - \int_{0}^t \tilde\bbE_{\tilde\rmP}\left[\nabla^2 W(X_s-\tilde X_s)\cdot(\Psi_s(X_s) - \Psi_s(\tilde{X}_s) ) \right]\, \De s
\ee
is a continuous martingale under $\rmP$ on $[0,T[$ and satisfies $\bbE_{\rmP} \left[ \int_{0}^T |M_t|^2 \De t \right]<+\infty$.
\end{prop}

To carry out the proof, we will use a well-known characterization of martingales (see e.g.\ \cite{emery1988cherchant}) which is as follows: an adapted process $(M_t)_{t\in[0,T]}$ such that $\bbE_{\rmP} \left[ \int_{0}^T |M_t|^2 \De t \right]<+\infty$ is a martingale in $[0,T[$ under $\rmP$ if and only if
\be\label{martortloop} \textstyle 
\bbE_{\rmP}\left[ \int_{0}^T M_t h_t \De t \right]=0
\ee
for all adapted processes $(h_t)_{t\in[0,T]}$  such that 
\be\label{loopcond} \textstyle 
\bbE_{\rmP} \left[ \int_{0}^T |h_t|^2 \De t \right]<+\infty, \quad \text{and} \,   \int_{0}^T h_t \,\De t=0 \quad \rmP-\text{a.s.}
\ee

\begin{proof}
Define $(M_t)_{t\in[0,T]}$ via \eqref{martingalecorrection1}. Using \eqref{boundedhess},\eqref{finitent2} and \eqref{secondmoment} we get  that  $\bbE_{\rmP}[\int_{0}^T|M_t|^2\De t] <+\infty$. Therefore, using the characterization of martingales \cite[pg.\ 148-149]{emery1988cherchant} in order to show that $M_t$ is a martingale on $[0,T[$ we need to show \eqref{martortloop} for all adapted processes $(h_t)_{t\in[0,T]}$ satisfying \eqref{loopcond}. By a standard density argument, one can show that it suffices to obtain \eqref{martortloop} under the additional assumption that $(h_t)_{t\in[0,T]}$ is bounded and Lipschitz, i.e.
\be\label{Lipschitzperturb} \textstyle 
\forall t\in [0,T], \,  \omega,\bar{\omega}\in\Omega, \quad \sup_{s\in[0,t]} |h_s(\omega)-h_s(\bar{\omega})| \leq C \sup_{s\in[0,t]} |\omega_s-\bar{\omega_s}|,\quad \sup_{t\in [0,T]}|h_t(\omega)|  \leq C,
\ee
for some $C>0$. Consider now a process $(h_t)_{t\in[0,T]}$ satisfying \eqref{loopcond} and \eqref{Lipschitzperturb} and for $\varepsilon>0$ define the shift transformation 
\be\label{shiftdef} \textstyle 
\tau^{\varepsilon}:\Omega \longrightarrow \Omega, \quad \tau^{\varepsilon}_t(\omega) = \omega_t +\varepsilon \int_{0}^t h_s(\omega) \De s.
\ee
Under the current assumptions, $\tau^{\varepsilon}$ admits an adapted inverse $Y^{\varepsilon}$, i.e.\ there exists an adapted process $(Y^{\varepsilon}_{t})_{t\in[0,T]}$ such that 
\be\label{shiftinvert2} \textstyle 
 \rmP-\text{a.s.}\quad \tau^{\varepsilon}_t (Y^{\varepsilon} (\omega))=Y^{\varepsilon}_t (\tau^{\varepsilon} (\omega)) = \omega_t \quad \forall t\in [0,T].
\ee
 {Indeed, since $\cH(\rmP|\Gamma(\rmP))<+\infty$, Lemma \ref{entropyequivalence} in the appendix yields that $\rmP\ll\rmR^{\mui}$; this entitles us to apply Lemma \ref{shiftinvert} in the same section, providing the existence of the inverse $Y^{\varepsilon}$.}

If we set $\rmP^{\varepsilon}=\rmP\circ(\tau^{\varepsilon})^{-1}$ we have that $\rmP^{\varepsilon}\in \cP_{1}(\Omega)$ is admissible for \eqref{nonlinearSP}, thanks to \eqref{loopcond}. Moreover, Lemma \ref{correctionexistence intro} and \eqref{shiftinvert2} imply that 
\bes  \textstyle 
X_t - \int_{0}^t \Big( \varepsilon h_s(Y^{\varepsilon})+\Psi_s(Y^{\varepsilon}_s)-\nabla W \ast \rmP_s(Y^{\varepsilon}_s)  \Big) \De s
\ees
is a Brownian motion under $\rmP^{\varepsilon}$. Combining \eqref{finitent2}, \eqref{Lipschitzperturb}
and \eqref{boundedhess} we get that
\be\label{martingalecorrection4} \textstyle 
\frac{1}{2}\bbE_{\rmP^{\varepsilon}}\left[ \int_{0}^T\big|\Psi_t(Y^{\varepsilon}_t) +\varepsilon h_t(Y^{\varepsilon}) -\nabla W \ast \rmP_t(Y^{\varepsilon}_t)+\nabla W\ast \rmP^{\varepsilon}_t(X_t) \big|^2\De t\right]
 <+\infty. 
\ee
Lemma \ref{driftregularity} grants that $\bar b(t,x)=-\nabla W \ast \rmP^{\varepsilon}_t(x)$ fulfills the hypothesis of Lemma \ref{finitent} and \eqref{martingalecorrection4} allows to use the implication $(ii)\Rightarrow (i)$ which yields that $\cH(\rmP^{\varepsilon}|\Gamma(\rmP^{\varepsilon}))$ is finite and equals the left hand side of \eqref{martingalecorrection4}. Using the definition of $\rmP^{\varepsilon}$, we 
can rewrite $\cH(\rmP^{\varepsilon}|\Gamma(\rmP^{\varepsilon}))$ as 
\bes \textstyle 
\frac{1}{2}\bbE_{\rmP}\left[ \int_{0}^T| \varepsilon h_t+\Psi_t(X_t) +\nabla W\ast \rmP^{\varepsilon}_t(\tau^{\varepsilon}_t)-\nabla W \ast \rmP_t(X_t) |^2\De s\right].
\ees
Imposing optimality of $\rmP$ and letting $\varepsilon$ to zero, using Taylor's expansion
\beas
0 &\leq &  \textstyle \liminf_{\varepsilon\rightarrow 0} \frac{\cH(\rmP^{\varepsilon}|\Gamma(\rmP^{\varepsilon}))-\cH(\rmP|\Gamma(\rmP)) } {\varepsilon}\\
&=&  \textstyle \bbE_{\rmP} \left[ \int_{0}^T\Psi_t(X_t) \cdot \left(  h_t + \tilde\bbE_{\tilde{\rmP}} \left[ \nabla^2 W(X_t-\tilde{X}_t) \cdot \int_{0}^t h_s-\tilde{h}_s \De s  \right]         \right)  \De t\right].
\eeas
{In the above equation, $(\tilde{X}_t,\tilde{h}_t)_{t\in[0,T]}$ is an independent copy of $(X_t,h_t)_{t\in[0,T]}$ defined on some probability space $(\tilde{\Omega},\tilde{\cF},\tilde{\rmP})$ and $\tilde{\bbE}_{\tilde{\rmP}}$ denotes the expectation on $(\tilde{\Omega},\tilde{\cF},\tilde{\rmP})$. Moreover, the exchange of limit and expectation is justified by \eqref{lipschitz}, \eqref{finitent2} and the dominated convergence theorem.} Using the symmetry of $W$, and taking $\pm h$, we can rewrite the latter condition as
$$0= \textstyle \bbE_{\rmP}\left[ \int_{0}^T\Psi_t(X_t) \cdot   h_t \De t \right] + \bbE_{\rmP} \left[ \int_{0}^T\tilde\bbE_{\tilde{\rmP}}\left[\ (\Psi_t(X_t)-\Psi_t(\tilde{X}_t)) \cdot  \nabla^2 W(X_t-\tilde{X}_t)\right] \cdot \int_{0}^t h_s\De s \,   \De t\right].$$
By integration by parts and the boundary condition \eqref{loopcond}
, we arrive at 
\bes  \textstyle 
0= \bbE_{\rmP} \left[ \int_{0}^T \left(  \Psi_t(X_t) - \int_{0}^t\tilde\bbE_{\tilde{\rmP}}\left[(\Psi_s(X_s)-\Psi_s(\tilde{X}_s)) \cdot   \nabla^2 W(X_s-\tilde{X}_s)  \right]   \De s \right) \cdot h_t           \De t \right],
\ees
proving the desired martingale property. {By \cite[Theorem IV.36.5]{RW00} we know that a martingale in an augmented Brownian filtration admits a continuous version. { Using again Lemma \ref{entropyequivalence} we have that $\rmP\ll \rmR^\mu$}, and we so obtain a continuous version of our martingale \eqref{martingalecorrection1}, and a fortiori of $\Psi_t(X_t)$.}
\end{proof}

\subsection{Benamou-Brenier formulation}\label{sec ben bren proofs}
We finally turn to the Benamou-Brenier formulation. Recall that $\nonlinentcost(\mui,\muf)$ denotes the optimal value of the mean field Schr\"odinger problem. 
We define the set $\mathcal A$ as the collection of all absolutely continuous curves $(\mu_t)_{t\in[0,T]}\subset \cP_2(\RD)$ (see Section \ref{OTsection}) such that 
\begin{align*}
(t,z)\mapsto\nabla \log \mu_t(z) &\in L^2(\De\mu_t\De t),\\
(t,z)\mapsto\nabla W \ast \mu_t(z) &\in L^2(\De\mu_t\De t).
\end{align*}
Recall from the introduction the problem
\be\label{eq cost BB}
\nonlinentcost^{BB}(\mui,\muf):=\inf_{\substack{(\mu_t)_{t\in[0,T]}\in\mathcal A, \\\partial_t\mu_t + \nabla \cdot(w_t\mu_t)=0 }} \, \frac{1}{2}\int\int \left |w_t(z)+  \frac{1}{2}\nabla \log \mu_t(z) + \nabla W \ast \mu_t(z)\right|^2\mu_t(\De z) \De t  
\ee
In \eqref{eq cost BB}, solutions to the continuity equation $\partial_t\mu_t + \nabla \cdot(w_t\mu_t)=0$ are meant in the weak sense. 

%
%

\begin{proof}[Proof of Theorem \ref{Teo BB intro}]
We first show that $\nonlinentcost(\mui,\muf) \geq \nonlinentcost^{BB}(\mui,\muf)$. To this end, we may assume that the l.h.s.\ if finite and denote $\rmP$ an optimizer. As established in Theorem \ref{martingalecorrection intro}, the drift of $X$ under $\rmP$ is equal to $$ \textstyle \int_0^t\Psi_s(X_s)-\nabla W\ast \rmP_s(X_s)\De s ,$$
where $\Psi\in \rmH_{-1}((\rmP_t)_{t\in[0,T]})$ and
$$ \textstyle \nonlinentcost(\mui,\muf) =  \frac{1}{2} \int\int |\Psi_t(z)|^2\rmP_t(\De z)\De t .$$
As we will see in Lemma \ref{OTlemma} and Remark \ref{Himpliestang}, the flow of marginals $(\rmP_t)_{t\in[0,T]}$ is absolutely continuous and its tangent velocity field $v$ is given by
$$ \textstyle v_t(z):= - \nabla W\ast \rmP_t(z)+ \Psi_t(z) -\frac{1}{2} \nabla\log \rmP_t.$$
Hence
\begin{align*}\nonlinentcost(\mui,\muf) &=  \textstyle \frac{1}{2} \int\int \left|v_t(z)+ \frac{1}{2} \nabla\log \rmP_t+\nabla W\ast \rmP_t(z)\right |^2\rmP_t(\De z)\De t .
\end{align*}
We conclude the desired inequality by noticing that $ \nabla \log \rmP_t \in L^2(\De\rmP_t\De t)$ and $\nabla W \ast \rmP_t \in L^2(\De\rmP_t\De t)$. To wit, the first statement follows from \cite[Thm 3.10]{follmer1986time} combined with Lemma \ref{entropyequivalence} in our appendix, and the second from \eqref{secondmoment} used with $\bar{b}=-\nabla W \ast \rmP_t(z)$.
We now establish $\nonlinentcost(\mui,\muf) \leq \nonlinentcost^{BB}(\mui,\muf)$, so we may assume that $(\mu_t)_{t\in[0,T]}$ is feasible for the r.h.s.\ and leads to a finite value. Denote by $\tilde v$ its tangent velocity field. We define $\Phi_t(z):= \tilde v_t(z)+\frac{1}{2}\nabla \log\mu_t(z)$, so from the continuity equation for $(\mu_t)_{t\in[0,T]}$ we deduce the following equation in the distributional sense
$$ \textstyle \partial_t \mu_t+\nabla\cdot(\mu_t \Phi_t)-\frac{1}{2}\Delta \mu_t=0.$$
{Observing} that $\Phi\in \rmH_{-1}((\mu_t)_{t\in[0,T]})$, we may apply the equivalence ``(a) iff (c)'' in \cite[Theorem 3.4]{cattiaux1995large} \footnote{That is for the construction of a Nelson process with marginals $(\mu_t)_{t\in[0,T]}$, with respect to the reference measure given by Wiener started at $\mui$.}. We thus obtain a measure $\rmP$ whose marginals are exactly $(\mu_t)_{t\in[0,T]}$, and by the uniqueness statement in \cite[Theorem 3.4]{cattiaux1995large} we also know that the drift of $X$ under $\rmP$ is precisely $ \Phi_s(X_s) $. Hence
\begin{align*} \textstyle 
\frac{1}{2}\int\int \left |\tilde v_t(z)+  \frac{1}{2}\nabla \log \mu_t(z) + \nabla W \ast \mu_t(z)\right|^2\mu_t(\De z) \De t &=  \textstyle \frac12\int\int \left |\Phi_t(z) + \nabla W \ast \mu_t(z)\right|^2\mu_t(\De z) \De t \\
&= \textstyle \frac{1}{2}\bbE_{\rmP}\left[\int_0^T|\Phi(X_t)+\nabla W\ast\rmP_t(X_t)|^2 \De t\right]\\
&\geq \nonlinentcost(\mui,\muf), 
\end{align*}
 where the inequality follows from the equivalent expression of  $\nonlinentcost(\mui,\muf)$ given in \eqref{equivSch}.
 
 We have proven $\nonlinentcost(\mui,\muf) = \nonlinentcost^{BB}(\mui,\muf)$, and the other statements follow from the previous arguments.
\end{proof}

\subsection{Schr\"odinger potentials and mean field PDE system: proofs}\label{sec schr pot}

\begin{proof}[Proof of Corollary \ref{mean field PDE sys}]
We know by Theorem \ref{correctionexistence intro} that $\Psi$ belongs to $\rmH((\mu_t))_{t\in[0,T]}$. The regularity hypothesis imposed on $\Psi_t(x)$ and $\mu_t(x)$ allow us to conclude that $\Psi$ is a true gradient, i.e. there exist $\psi$ such that $\Psi_t(x)=\nabla \psi_t(x)$ for all $(t,x)\in [0,T]\times \RD$. Lemma \ref{correctionexistence intro} together with Theorem \ref{martingalecorrection intro} yield that $\mu_t$ is a weak solution of the Fokker Planck equation in \eqref{FBSDE master}. Because of the regularity assumptions we made on $\Psi$ and $\mu$, we can conclude that $\mu_t$ is indeed a  classical solution. For the same reasons, we can turn the martingale condition \eqref{martingalecorrection intro} into the system of PDEs
\bes  \textstyle 
\forall i=1,\ldots,d \quad \partial_t \partial_{x_i}\psi_t(x) + \mathcal{L} (\partial_{x_i} \psi_t(x)) - \int_{\RD} \partial_{x_i}((\nabla W(x-\tilde{x}))\cdot ( \nabla\psi(x) -\nabla \psi(\tilde{x})) \mu_t(\De \tilde{x}) =0,
\ees
where $\mathcal{L}$ is the generator $\frac{1}{2}\Delta + (\nabla(-W \ast \mu_t + \psi_t))\cdot \nabla $. After some tedious but standard calculations we can rewrite the above as
\bes  \textstyle 
\partial_{x_i} \left(  \partial_t \psi_{t}(x) + \frac{1}{2}\Delta \psi_{t}(x) + \frac{1}{2} |\nabla\psi_{t}(x)|^2 + \int_{\RD} \nabla W(x-\tilde{x}) \cdot( \nabla \psi_t(x)  -\nabla \psi_t(\tilde{x})\mu_t(\De \tilde{x}) ) \right) =0
\ees

Since $\psi$ is defined up to the addition of a function that depends on time only, the conclusion follows. 

\end{proof}

Corollary \ref{cor Schr pot} can be proven with a direct calculation using the definition of $\varphi_t$ and \eqref{FBSDE master}.
\section{Convergence to equilibrium and functional inequalities: proofs}\label{quant res proofs}

In this part we complement the discussion undertaken in Section \ref{sec exp eq functional ineq} and provide proofs for the results stated therein. This section is organized as follows:
\bei 
\item Subsections \ref{sub corrector bound},\ref{OTsection},\ref{reversalsection} are devoted to stating and proving some preparatory results that we shall use at different times in the proofs of the main results.
\item In Subsection \ref{sub fun eq proofs} we prove Theorem \ref{consquantity1intro} and Theorem \ref{nonlinearentropybound1intro} together with its corollaries, ie.\ the Talagrand inequality (Corollary \ref{talagrand}) and the HWI inequality (Corollary \ref{HWI}). 
\item Finally, in Subsection \ref{con eq proofs} we prove Theorem \ref{exp conv eq} and Theorem \ref{thm exp mkv intro}.
\eei

In all the lemmas and theorems in this subsection we always assume \eqref{boundedhess}-\eqref{marginalhyp} to hold, and throughout $\rmP,\alpha^{\rmP},\Psi,M$ are as given in Theorem \ref{martingalecorrection intro}. We refer to Sections \ref{sec mfs intro} and \ref{sec exp eq functional ineq} for any unexplained notation.

\subsection{Exponential upper bound for the corrector}\label{sub corrector bound}
Recall that we called $\Psi$ the \emph{corrector}. The goal of this part is to quantify the size of the corrector, as stated in Lemma \ref{correctionboundlemma} below. Before doing this we prove two preliminary lemmas. As usual, we denote by $\langle \cdot \rangle$ the quadratic variation of a semimartingale.

\begin{lemma}\label{L2bounds}
We have
\be\label{boundedL2martingale}
\forall t\in[0,T[, \quad \bbE_{\rmP}[|M_t|^2]=\bbE_{\rmP}[\langle M\rangle_t]<+\infty.
\ee
Moreover the function $t\mapsto \bbE\left[  \langle M\rangle_t \right]$ is continuous on $[0,T[$ and
\be\label{boundedL2correction} \textstyle 
\forall t\in[0,T[, \quad \sup_{s\in[0,t]}\, \bbE_{\rmP}[|\Psi_s(X_s)|^2]<+\infty
\ee
\end{lemma}
\begin{proof}
We have shown at Theorem \ref{martingalecorrection} that $\bbE_{\rmP}\left[ \int_0^T|M_t|^2 \De t \right]<+\infty$ which gives that $\bbE_{\rmP}\left[ |M_t|^2 \right]<+\infty$ for almost every $t\in[0,T[$. But since $\bbE_{\rmP}\left[ |M_t|^2 \right]$ is an increasing function of $t$, we get $\bbE_{\rmP}\left[ |M_t|^2\right]<+\infty$ for all $t\in[0,T[$. To complete the proof of \eqref{boundedL2martingale} it suffices to observe that by definition of quadratic variation and since $M_t$ is an $L^2$-martingale on $[0,T[$, we have $\bbE_{\rmP}\left[|M_t|^2 \right]=\bbE_{\rmP}[\langle M_t\rangle]$. To prove the continuity of $t\mapsto\bbE_{\rmP}[\langle M\rangle_t]$ we start by observing that since $M_t$ is a continuous martingale, then $\langle M\rangle_t$ has continuous and increasing paths. Thus, we obtain by monotone convergence that  $\bbE_{\rmP}[\langle M\rangle_{t+h}] \rightarrow \bbE_{\rmP}[\langle M\rangle_{t}] $ as $h\downarrow 0$, which gives the desired result. The proof of \eqref{boundedL2correction} follows from \eqref{martingalecorrection1}, the bounded Hessian of $W$ (see\eqref{boundedhess}) and the first part of Theorem \ref{martingalecorrection intro}.
\end{proof}

\begin{lemma}\label{linearbar}
The function $t\mapsto \bbE_{\rmP}[X_t]$ is linear, the function $t\mapsto \bbE_{\rmP}[\Psi_t(X_t)]$ is constant, and
\bes
\forall t\in[0,T[,\quad \bbE_{\rmP}[X_t]=\bbE_{\rmP}[X_0] + \bbE_{\rmP}[\Psi_0(X_0)] t 
\ees
\end{lemma}

\begin{proof}
Using the symmetry of $W$ and the martingale property \eqref{martingalecorrection1} it is easily derived that $\bbE_{\rmP}[\Psi_t(X_t)]$ is constant as a function of $t$. Therefore we get for all $t\in[0,T]$
\bes
\bbE_{\rmP}[X_t]=\bbE_{\rmP}[X_0]-\int_{0}^t \bbE_{\rmP}[\nabla W\ast \rmP_s(X_s)]\De s + \bbE_{\rmP}[\Psi_0(X_0)] t 
\ees
Using again the symmetry of $W$ we get that $\int_{0}^t \bbE_{\rmP}[\nabla W\ast \rmP_s(X_s)]\De s =0$, from which the conclusion follows.
\end{proof}

We can now provide some key estimates on the corrector:



\begin{lemma}\label{correctionboundlemma}
Assume \eqref{boundedhess}-\eqref{same mean ass}. If $\rmP$ is an optimizer for \eqref{nonlinearSP} and $\Psi$ the associated corrector,  then for any $t\in(0,T)$ we have
\be\label{correctionbound1} \textstyle 
\frac{1}{2}\bbE_{\rmP}\left[\int_{0}^t |\Psi_s(X_s)|^2\De s\right]\leq  \frac{\exp(2\curv t)-1}{\exp(2\curv T)-1} \nonlinentcost(\mu^{\mathrm{in}},\mu^{\mathrm{fin}}),
\ee
and
\be\label{correctionbound3} \textstyle 
\frac{1}{2}\bbE_{\rmP}\left[ |\Psi_t(X_t)|^2 \right] \leq  \frac{2\curv\, \nonlinentcost(\mu^{\mathrm{in}},\mu^{\mathrm{fin}})}{\exp(2\curv(T-t))-1}.
\ee

\end{lemma}

\begin{proof}
Consider the function $t\mapsto \varphi(t)$ defined by
\bes  \textstyle 
\varphi(t) = \frac12\bbE_{\rmP} \left[ \int_0^t |\Psi_s(X_s)|^2 \De s \right].
\ees
Fubini's Theorem allows to interchange the time integral and the expectation to get that $\varphi$ is an absolutely continuous function with derivative
\be\label{correctionbound5}
\varphi'(t)=\bbE_{\rmP} \left[|\Psi_t(X_t)|^2  \right].
\ee
From It\^o's formula and Theorem \ref{martingalecorrection intro} we get that for all $t\in[0,T[$
\begin{align*} \textstyle 
|\Psi_t(X_t)|^2-|\Psi_0(X_0)|^2 = & \textstyle  2\int_0^t \Psi_r(X_r) \cdot \De M_r + \\ &  \textstyle 2\int_0^t\Psi_r(X_r)\cdot \tilde{\bbE}_{\tilde{\rmP}}[\nabla^2 W(X_r-\tilde{X}_r)\cdot (\Psi_r(X_r)-\Psi_r(\tilde{X_r}))] \, \De r + \langle M \rangle_t.
\end{align*}
We observe that the fact that $M_t$ is a martingale together with \eqref{boundedL2correction} and \eqref{boundedL2martingale} make sure that $\bbE_{\rmP}\left[ \int_0^t \Psi_r(X_r)\cdot \De M_r \right]=0$. Thus, taking expectation on both sides of the above equation yields
\be\label{correctionvscost1} \textstyle 
\varphi'(t)-\varphi'(0) = \bbE_{\rmP}\left[2 \int_0^t\Psi_r(X_r)\cdot \tilde{\bbE}_{\tilde{\rmP}}[\nabla^2 W(X_r-\tilde{X}_r)\cdot (\Psi_r(X_r)-\Psi_r(\tilde{X_r}))] \, \De r\right] + \bbE_{\rmP}[\langle M \rangle_t].
\ee
Because of \eqref{boundedL2correction} we can use Fubini's Theorem and write
\begin{align}\label{correctionvscost2}\nonumber  &  \textstyle \bbE_{\rmP}\left[2 \int_0^t\Psi_r(X_r)\cdot \tilde{\bbE}_{\tilde{\rmP}}[\nabla^2 W(X_r-\tilde{X}_r)\cdot (\Psi_r(X_r)-\Psi_r(\tilde{X_r}))] \, \De r\right]  \\
\nonumber = &2  \textstyle \int_0^t \bbE_{\rmP}\left[ \Psi_r(X_r)\cdot \tilde{\bbE}_{\tilde{\rmP}}[\nabla^2 W(X_r-\tilde{X}_r)\cdot (\Psi_r(X_r)-\Psi_r(\tilde{X_r}))]\, \right]\De r\\
= &  \textstyle \int_0^t \bbE_{\rmP \otimes \tilde{\rmP}}\left[ (\Psi_r(X_r)-\Psi_r(\tilde{X}_r)) \cdot\nabla^2 W(X_r-\tilde{X}_r)\cdot (\Psi_r(X_r)-\Psi_r(\tilde{X}_r)) \, \right]\De r,
\end{align} 
where we used the symmetry of $W$ to obtain the last expression. Plugging it back in
\eqref{correctionvscost1} and using that  $t\mapsto\varphi'(t)$ is
\bei 
\item continuous on $[0,T[$ because so are  \eqref{correctionvscost2} and $\bbE[\langle M \rangle_t]$ (cf.\ Lemma \ref{L2bounds}),
\item increasing on $[0,T[$ since $W$ is convex and the quadratic variation is an increasing process,
\eei
we conclude that $t\mapsto\varphi'(t)$ is absolutely continuous on the same interval. Moreover, using the $\curv$-convexity of $W$ and again the fact that the quadratic variation is an increasing process we get 
\be\label{curvaturebound}
\varphi''(t) \geq 2 \curv\, \bbE_{\rmP}[|\Psi_t(X_t)|^2 ]=2\curv \, \varphi'(t)
\ee
where to establish the last inequality we used that the hypothesis on $\mu^{\mathrm{in}}$ and $\mu^{\mathrm{fin}}$ together with Lemma \ref{linearbar} imply $\bbE_{\rmP}[\Psi_t(X_t)]=0$. 
The bound \eqref{correctionbound1} follows by integrating the differential inequality \eqref{curvaturebound} as done for instance in Lemma \ref{differentialinequality} in the Appendix, and observing that 
$\frac{1}{2}\bbE_{\rmP}\left[\int_{0}^T|\Psi_r(X_r)|^2\De r \right]=\cH(\rmP|\Gamma(\rmP))$. To prove \eqref{correctionbound3}, we begin by observing that \eqref{curvaturebound} also yields that 
\be\label{correctionbound4}
\forall s\in [t,T], \quad \bbE_{\rmP}\left[|\Psi_s(X_s)|^2 \right]\geq \exp(2\curv(s-t)) \bbE_{\rmP}\left[|\Psi_t(X_t)|^2 \right].
\ee
Next, by definition of entropic cost we get the trivial bound
\bes  \textstyle 
\nonlinentcost(\mu^{\mathrm{in}},\mu^{\mathrm{fin}}) =\frac{1}{2}\bbE_{\rmP}\left[ \int_{0}^T |\Psi_s(X_s)|^2\De s \right] \geq \frac{1}{2}\bbE_{\rmP}\left[ \int_{t}^T |\Psi_s(X_s)|^2\De t \right]
\ees
The desired conclusion follows by plugging \eqref{correctionbound4} in the above equation and some standard calculations.
\end{proof}



\subsection{First derivative of $\nent$}\label{OTsection}
We compute the first derivative of $\nent$ along the marginal flow of $\rmQ$, assuming that $\cH(\rmQ|\Gamma(\rmQ))<+\infty$ and that $\rmQ$ is Markovian. To do this, we use an approach based on optimal transport, and some results of \cite{follmer1986time}. To be self-contained, we recall the basic notions of optimal transport we need to state the results. We refer to \cite{ambrosio2008gradient} for more details.

\paragraph{Tangent space} Let $\mu\in\cP_2(\RD)$. The tangent space $\tang_{\mu}\cP_2$ at $\mu$ is the closure in $L^2_{\mu}$ of 
\bes
\left\{ \nabla \psi; \psi \in \cC^{\infty}_c(\RD) \right\}.
\ees
Since $L^2_{\mu}$ is an Hilbert space, given an arbitrary $\Psi\in L^2_{\mu}$, there exists a unique projection $\Pi_{\mu}(\Psi) $ of $\Psi$ onto $\tang_{\mu}\cP_2(\RD)$. 

\paragraph{Absolutely continuous curves and velocity field}

Following \cite[Th 8.3.1]{ambrosio2008gradient}, we say that a curve $(\mu_t)_{t\in[0,T]} \subseteq \cP_2(\RD)$ is \emph{absolutely continuous} if there exists a Borel measurable vector field $(t,z)\mapsto w_t(z)$ such that 
\bei 
\item $(w_t)_{t\in[0,T]}$ solves (in the sense of distributions) the continuity equation
\be\label{conteq}
\partial_t \mu_t + \nabla \cdot (w_t \mu_t )=0.
\ee

\item $w_t$ satisfies the integrability condition
\bes  \textstyle 
\int_{0}^T\left(\int_{\RD}|w_t(z)|^2\mu_t(\De z)\right)^{1/2} \De t<+\infty.
\ees
\eei

Consider an absolutely continuous curve $(\mu_t)_{t\in[0,T]}$. It is a consequence of the results in Chapter 8, and in particular of Proposition 8.4.5 of \cite{ambrosio2008gradient}, that there exist a unique Borel measurable vector field $v_t(z)$ solving \eqref{conteq} and such that $z\mapsto v_t(z)$ belongs to the tangent space $\tang_{\mu_t}\cP_2$ for almost every $t\in[0,T]$. We call such $v_t$ the \emph{(tangent) velocity field of $(\mu_t)_{t\in[0,T]}$}. 

\begin{remark}\label{Himpliestang}
Let $(\mu_t)_{t\in[0,T]}$ be an absolutely continuous curve and $w_t(z)$ be in $\rmH_{-1}((\mu_t)_{t\in[0,T]})$. It is rather easy to see that $z\mapsto w_t(z)$ belongs to $\tang_{\mu_t}\cP_2$ for almost every $t\in[0,T]$.
\end{remark}

Throughout the rest of the paper, if $\rmQ\in \cP(\Omega)$ is such that $\cH(\rmQ | \Gamma(\rmQ))<+\infty$, we say that $\rmQ$ is \emph{Markovian} if $\alpha^{\rmQ}_t$ is $\sigma(X_t)$-measurable for all $t\in[0,T]$, $(\alpha^{\rmQ}_t)_{t\in[0,T]}$ being defined by \eqref{finitent2}. In that case we write $\Xi^{\rmQ}_t(X_t)$ instead of $\alpha^{\rmQ}_t$.

\begin{lemma}\label{OTlemma}
Let $ \rmQ $ be such that $\cH(\rmQ|\Gamma(\rmQ))<+\infty$ and Markovian. Then

\bei
\item[(i)] $(\rmQ_t)_{t\in[0,T]}$ is an absolutely continuous curve. Its tangent velocity field is given by
\be\label{velfield1} \textstyle 
v_t(z) = -\nabla W \ast \rmQ_t(z) + \Pi_{\rmQ_t}(\Xi^{\rmQ}_t(z)) - \frac{1}{2} \nabla \log \rmQ_t(z).
\ee
Moreover, 
\be\label{finitenrg} \textstyle 
\int_{0}^T \int_{\RD} |v_t|^2 \De \rmQ_t \De t <+\infty.
\ee
\item[(ii)] The function $t\mapsto \nent(\rmP_t)$ is absolutely continuous and
\be\label{chainrule} \textstyle 
\forall 0\leq s\leq t, \quad \nent(\rmQ_t)-\nent(\rmQ_s) = \int_s^t \int_{\RD}\big(\nabla \log \rmQ_r + 2 \nabla W\ast \rmQ_r \big)(z) \cdot v_r(z)  \rmQ_r(\De z) \,\De r.
\ee
\eei
\end{lemma}

\begin{proof}$ $

\noindent 
\underline{Proof of (i)}  To show that $(\rmQ_t)_{t\in[0,T]}$ is absolutely continuous it suffices to show that there exists a distributional solution of the continuity equation 
\be\label{coneq}
\partial_t \rmQ_t + \nabla \cdot (w_t \rmQ_t)=0
\ee
with the property that 
\be\label{coneq3} \textstyle 
\int_{0}^T \left(\int_{\RD} |w_t(z)|^2 \rmQ_t(z)\right)^{1/2} \De t <+\infty.
\ee
Let now $\varphi\in \cC^{\infty}_c(]0,T[ \times \RD)$. Using It\^o's formula and taking expectation we obtain
\be\label{coneq2} \textstyle 
\int_{0}^T \int_{\RD}\Big(\nabla \varphi(t,z) \big( -\nabla W \ast \rmQ_t(z) + \Xi^{\rmQ}_t(z) \big) + \frac{1}{2}\Delta \varphi(t,z) + \partial_t \varphi(t,z) \Big) \, \rmQ_t(\De z)=0 .
\ee
Lemma \ref{entropyequivalence} in the appendix grants that under the current assumptions $\cH(\rmQ|\rmR^{\mui})<+\infty$, where $\rmR^{\mui}$ is the Wiener measure started at $\mu^{\mathrm{in}}$. But then, using \cite[Thm 3.10]{follmer1986time}\footnote{Strictly speaking, F\"ollmer's result is only concerned with the case $\mui=\delta_0$. However, a simple adaptation of his argument show that its validity extends to any $\mui$ satisfying \eqref{marginalhyp}.} we obtain that
$\log \rmQ_t$ is an absolutely continuous function for almost every $t$ and that $(t,z)\mapsto \nabla \log \rmQ_t(z)$ belongs to $\rmH_{-1}((\rmQ_t)_{t\in[0,T]})$. Therefore we can use integration by parts in \eqref{coneq2} to obtain
\bes \textstyle 
\forall t\in[0,T], \quad \frac{1}{2}\int_{\RD} \nabla \varphi(t,z) \rmQ_t(\De z) = -\frac{1}{2}\int_{\RD} \Delta \log \varphi(t,z) \cdot \nabla \log \rmQ_t(z) \,\rmQ_t(\De z)
\ees
which gives, using the definition of the projection operator $\Pi_{\rmQ_t}$,  that the rhs of \eqref{velfield1} solves the continuity equation in the sense of distributions. Next, we observe that \eqref{finitent2} grants that $\Pi_{\rmQ_t}(\Xi^{\rmQ}(t,z))\in\rmH_{-1}((\rmQ)_{t\in[0,T]})$.  We have already shown that $\nabla \log \rmQ_t \in \rmH_{-1}((\rmQ)_{t\in[0,T]})$, and \eqref{secondmoment} used with $\bar{b}=-\nabla W \ast \rmQ_t(z)$ yields that $-\nabla W\ast \rmQ_t(z)\in\rmH_{-1}((\rmQ)_{t\in[0,T]})$. Thus $v_t(z)\in \rmH_{-1}((\rmQ)_{t\in[0,T]})$ as well, which gives \eqref{coneq3} and \eqref{finitenrg}. Finally, Remark \ref{Himpliestang} yields that $(v_t)_{t\in[0,T]}$ is indeed the tangent velocity field.

\noindent\underline{Proof of (ii)}
From point (i) we know that $z\mapsto \nabla \log \rmQ_t(z)$ belongs to $L^2_{\rmQ_t}$ for almost every $t$; this implies that $\nabla \log \rmQ_t + 2 \nabla W \ast \rmQ_t $ belongs to the subdifferential of $\nent$ at $\rmQ_t$ for almost every $t$ (see e.g. \cite[Thm. 10.4.13]{ambrosio2008gradient}). The chain rule \cite[sec. E, pg. 233-234]{ambrosio2008gradient} gives the desired result \eqref{chainrule}, provided its hypothesis are verified. 
 We have to check that (a) $(\rmQ_t)_{t\in[0,T]}$ is an absolutely continuous curve and $\nent(\rmQ_t)<+\infty$ for all $t\in[0,T]$, (b) $\nent(\cdot)$ is displacement $\lambda$-convex for some $\lambda \in \R$, and (c) that \bes  \textstyle 
\int_{0}^T \left(\int_{\RD}|v_t|^2 \De \rmQ_t\right)^{1/2}\left( \int_{\RD} \Big|\nabla \log \rmQ_t + 2 \nabla W\ast \rmQ_t\, \Big|^2 \De \rmQ_t \right)^{1/2} \De t <+\infty.
\ees
To wit, (a) follows from point (i) and the fact that $\cH(\rmQ|\Gamma(\rmQ))<+\infty$, and (b) is a consequence of displacement convexity of the entropy and \eqref{boundedhess}. Finally, (c) is granted by \eqref{finitenrg} and the fact that $\nabla \log \rmQ_t(z)+2 \nabla W\ast \rmQ_t(z)$ belongs to $\rmH_{-1}((\rmQ_t)_{t\in[0,T]})$ (see the proof of (i)).

\end{proof}

\subsection{Time reversal}\label{reversalsection}

For $\rmQ\in \cP(\Omega)$ the time reversal $\hat{\rmQ}$ is the law of the time reversed process $(X_{T-t})_{t\in[0,T]}$. In this section we derive an expression for $\cH(\hat{\rmQ}|\Gamma(\hat{\rmQ}))$ and use it to derive the bound \eqref{revcorrectionbound} below, which plays a fundamental role in the proof of Theorem \ref{nonlinearentropybound1intro}.

\begin{prop}\label{entropytimereversal3}
Let $\rmQ\in\cP_{1}(\Omega)$ be Markovian and such that $\cH(\rmQ|\Gamma(\rmQ))<+\infty$. 
\bei
\item[(i)] If $\rmQ_0=\mu^{\mathrm{in}}$, $\rmQ_T=\mu^{\mathrm{fin}}$ then $\cH(\hat{\rmQ}|\Gamma(\hat{\rmQ}))<+\infty$ as well and 
\be\label{entropytimereversal1}
\cH(\hat{\rmQ}|\Gamma(\hat{\rmQ})) =\cH(\rmQ|\Gamma(\rmQ))+\nent(\mu^{\mathrm{in}})-\nent(\mu^{\mathrm{fin}}) 
\ee
\item[(ii)] If $\rmQ_0=\mu^{\mathrm{fin}}$, $\rmQ_T=\mu^{\mathrm{in}}$ then $\cH(\hat{\rmQ}|\Gamma(\hat{\rmQ}))<+\infty$ as well and 
\be\label{entropytimereversal2}
 \cH(\rmQ|\Gamma(\rmQ))=\cH(\hat{\rmQ}|\Gamma(\hat{\rmQ}))+\nent(\mu^{\mathrm{in}})-\nent(\mu^{\mathrm{fin}}) 
\ee
\eei
\end{prop}

\begin{proof}
 We only prove (i), (ii) being completely analogous. Recalling (see Lemma \ref{entropyequivalence}) that $\cH(\rmQ|\Gamma(\rmQ))<+\infty$ implies $\cH(\rmQ|\rmR^{\mu^{\mathrm{in}}})<+\infty$, we can use \cite[Thm. 3.10, Eq. 3.9]{follmer1986time} to obtain that there exist a Borel measurable vector field $\hat{b}_t(x)$ such that 
\bes  \textstyle 
X_t - \int_{0}^t \hat{b}_s(X_s) \De s 
\ees
is a Brownian motion under $\hat{\rmP}$ and that
\be\label{timerevdrift}
 \hat{\rmQ}-\text{a.s.} \quad \hat{b}_t(X_t) = -b_{T-t}(X_t) + \nabla \log \rmQ_{T-t}(X_t)\quad \forall t\in[0,T],
\ee
where 
$b_t(z)$ is the drift of $\rmQ$, that, in view of Lemma \ref{correctionexistence intro} we write as $ -\nabla W \ast \rmQ_t  + \Xi^{\rmQ}_t(z)$. Thus, we deduce that under $\hat{\rmQ}$ we have that
\bes  \textstyle 
X_t- \int_{0}^t -\nabla W\ast \hat{\rmQ}_s(X_s)+\hat{\Xi}^{\rmQ}_s(X_s)\De s
\ees
is a Brownian motion, where 
\be\label{timerevcorrection4}
\hat{\rmQ}-\text{a.s.}\quad \hat{\Xi}^{\rmQ}_t(X_t)= -\Xi^{\rmQ}_{T-t}(X_t) + \nabla \log \rmQ_{T-t}(X_t)+2\nabla W\ast \rmQ_{T-t}(X_t) \quad \forall t\in[0,T].
\ee

In the proof of Lemma \ref{OTlemma}, it was shown that $(\nabla \log \rmQ_{t})_{t\in[0,T]}$, $(\nabla W \ast \rmQ_{\cdot})_{t\in[0,T]}$ and $(\Xi^{\rmQ}_t)_{t\in[0,T]}$ are all in $\rmH_{-1}((\rmQ_t)_{t\in[0,T]})$. This implies that $(\hat{\Xi}^{\rmQ}_t)_{t\in[0,T]}\in \mathrm{H}_{-1}((\hat{\rmQ}_t)_{t\in[0,T]})$ as well. But then using $(ii)\Rightarrow(i)$ in Lemma \ref{finitent} for the choice $\bar{b}(t,z)=-\nabla W\ast \hat{\rmQ}_t(z)$ we get that $\cH(\hat{\rmQ}|\Gamma(\hat{\rmQ}))<+\infty$ and 
\bes \textstyle 
\cH(\hat{\rmQ}|\Gamma(\hat{\rmQ})) = \frac{1}{2}\bbE_{\hat{\rmQ}} \left[ \int_0^T|\hat{\Xi}^{\rmQ}_t(X_t)|^2 \De t \right].
\ees
Using \eqref{timerevcorrection4} in the above equation we get
\beas
\cH(\hat{\rmQ}|\Gamma(\hat{\rmQ})) &=&  \textstyle 
 \frac{1}{2}\bbE_{\hat{\rmQ}} \left[ \int_0^T|\Xi^{\rmQ}_{T-t}(X_t) -\big( \nabla \log \rmQ_{T-t}(X_t)+2\nabla W\ast  \rmQ_{T-t}(X_t)\big)|^2 \De t \right]\\
&=& \textstyle \frac{1}{2} \bbE_{\rmQ} \left[ \int_0^T|\Xi^{\rmQ}_t(X_t)|^2 \De t\right]\\
&+& \textstyle  \frac{1}{2}\bbE_{\rmQ} \Big[ \int_0^T\big(\nabla \log \rmQ_{t}(X_t)+2\nabla W\ast  \rmQ_{t}(X_t) \big)\cdot \big(-2 \Xi^{\rmQ}_{t}(X_t) +  \nabla \log \rmQ_{t}(X_t) 
\\&+& \textstyle 2\nabla W\ast  \rmQ_{t}(X_t)\big)\De t \Big]\\
&\stackrel{\eqref{velfield1}}{=}&  \textstyle  \cH(\rmQ|\Gamma(\rmQ)) -\bbE_{\rmQ}\left[\int_{0}^T \big(\nabla \log \rmQ_{t}(X_t)+2\nabla W\ast  \rmQ_{t}(X_t) \big) \cdot v_t(X_t)\De t\right].
\eeas
The conclusion follows from point (ii) of Lemma \ref{OTlemma}.
\end{proof}

A consequence of Proposition \ref{entropytimereversal3} is that optimality is preserved under time reversal.
\begin{lemma}\label{timerevoptimality}
Let $\rmP$ be an optimizer for \eqref{nonlinearSP}. Then $\hat{\rmP}$ optimizes
\bea\label{timerevnonlinearSP}
\inf \left\{\cH(\rmQ | \Gamma(\rmQ) ) \,:\,  \rmQ \in \cP_{1}(\Omega),\,  \rmQ_0=\mu^{\mathrm{fin}},\, \rmQ_T=\mu^{\mathrm{in}} \right \}.
\eea
\end{lemma}

\begin{proof}
Let us observe that since \eqref{marginalhyp} makes no distinction between $\mu^{\mathrm{in}}$ and $\mu^{\mathrm{fin}}$, the problem \eqref{timerevnonlinearSP} admits at least an optimal solution by Proposition \ref{existence}. Applying Proposition \ref{Markovopt} inverting the roles of $\mu^{\mathrm{in}},\mu^{\mathrm{fin}}$ we get that the optimizers of \eqref{timerevnonlinearSP} are Markovian. So it suffices to show that for any Markovian $\rmQ$ admissible for \eqref{timerevnonlinearSP} we have $\cH(\rmQ|\Gamma(\rmQ))\geq\cH(\hat{\rmP}|\Gamma(\hat{\rmP}))$. Take any such $\rmQ$. We have
\begin{align*}
\cH(\rmQ|\Gamma(\rmQ)) \stackrel{\text{Prop.} \ref{entropytimereversal3} (ii)}{=} \cH(\hat{\rmQ}|\Gamma(\hat{\rmQ}))
+ \nent(\mu^{\mathrm{in}})-\nent(\mu^{\mathrm{fin}})
&\stackrel{\text{Opt. of $\rmP$}}{\geq} \cH(\rmP|\Gamma(\rmP))
+  \nent(\mu^{\mathrm{in}})-\nent(\mu^{\mathrm{fin}})\\
&\stackrel{\text{Prop.} \ref{entropytimereversal3} (i)}{=}\cH(\hat{\rmP}|\Gamma(\hat{\rmP})).
\end{align*}
\end{proof}

\subsection{Functional inequalities: proofs}\label{sub fun eq proofs}

The goal of this section is to prove Theorem \ref{consquantity1intro} and Theorem \ref{nonlinearentropybound1intro} together with its corollaries. 

\subsubsection{Proof of Theorem \ref{nonlinearentropybound1intro}}
Using a time reversal argument, we prove the bound \eqref{revcorrectionbound} which is the key ingredient of the proof of Theorem \ref{nonlinearentropybound1intro} together with the bound for the correction term \eqref{correctionbound1}.

Tthe \emph{backward corrector} $\hat{\Psi}$ is obtained by the same argument used in Proposition \ref{entropytimereversal3}, replacing \eqref{nonlinearSP} with \eqref{timerevnonlinearSP} to obtain that there exist a Borel measurable vector field $\hat{\Psi}_t(z)\in \rmH_{-1}((\hat{\rmP}_t)_{t\in[0,T]})$ such that 
\bes \textstyle 
X_t -\int_{0}^t\left(-\nabla W\ast\hat{\rmP}_{s}(X_s) +\hat{\Psi}_s(X_s)\right) \De s
\ees
is a Brownian motion under $\hat{\rmP}$. Moreover, the following relation holds 
\be\label{timerevcorrection5}
\hat{\rmP}-\text{a.s.}\quad \hat{\Psi}_t(X_t)= -\Psi_{T-t}(X_t) + \nabla \log \rmP_{T-t}(X_t)+2\nabla W\ast \rmP_{T-t}(X_t) \quad \forall t\in[0,T].
\ee

\begin{lemma}\label{timerevcorrectionboundlemma}
Assume \eqref{W convexity ass},\eqref{same mean ass} and let $\rmP$ be an optimizer for \eqref{nonlinearSP}. Then
\bea\label{revcorrectionbound} \textstyle 
\nent(\rmP_r)+ \frac{1}{2}\bbE_{\rmP}\left[ \int_{r}^{T} |\Psi_s(X_s)|^2 \De s \right] \leq \frac{\exp(2\curv (T-r))-1}{\exp(2\curv T)-1} \nonlinentcost(\mui,\muf)  \\  \textstyle 
\nonumber + \frac{\exp(2\curv (T-r))-1}{\exp(2\curv T)-1}\nent(\mu^{\mathrm{in}}) +  \frac{\exp(2\curv T)-\exp(2\curv (T-r))}{\exp(2\curv T)-1} \nent(\mu^{\mathrm{fin}})
\eea
and 
\be\label{revcorrectionbound3} \textstyle 
\frac{1}{2}\bbE_{\rmP}\left[ |\hat{\Psi}_{T-r}(X_r)|^2 \right] \leq  \frac{2\curv\, \nonlinentcost(\mu^{\mathrm{fin}},\mu^{\mathrm{in}})}{\exp(2\curv r)-1}
\ee
hold for all $r\in(0,T)$.
\end{lemma}

\begin{proof}
 Using \eqref{velfield1} we can rewrite the above equation \eqref{timerevcorrection5} as

\be\label{timerevcorrection6}
\hat{\rmP}-\text{a.s.}\quad \hat{\Psi}_t(X_t)= \Psi_{T-t}(X_t) -2v_{T-t}(X_t) \quad \forall t\in[0,T]
\ee
From Proposition \ref{timerevoptimality} we also know that $\hat{\rmP}$ is optimal for \eqref{timerevnonlinearSP} and hence that $\cH(\hat{\rmP}|\Gamma(\hat{\rmP}))=\nonlinentcost(\muf,\mui)$. Therefore, by inverting again the roles of $\mui$ and $\muf$, we can use Lemma \ref{correctionboundlemma} for the problem \eqref{timerevnonlinearSP} setting $t=T-r$ to derive that
\be\label{revcorrectionboundlemma2} \textstyle 
\frac{1}{2}\bbE_{\hat{\rmP}}\left[ \int_{0}^{T-r} |\hat{\Psi}_s(X_s)|^2\De s \right] \leq \frac{\exp(2\curv (T-r))-1}{\exp(2\curv T)-1} \cH(\hat{\rmP}|\Gamma(\hat{\rmP})).
\ee
Thanks to \eqref{timerevcorrection6} we can write
\bea\label{revcorrectionboundlemma4} \textstyle 
\nonumber\frac{1}{2}\bbE_{\hat{\rmP}}\left[ \int_{0}^{T-r}|\hat{\Psi}_s(X_s)|^2\De s \right]=
\frac{1}{2} \bbE_{\rmP} \left[ \int_r^T|\Psi_s(X_s)|^2 \De s\right]-\bbE_{\rmP} \left[ \int_r^T\big(2\Psi_s(X_s)-2v_s(X_s)) v_s(X_s) \De s \right]\\  \textstyle 
\nonumber\stackrel{\eqref{velfield1}+\Psi \in \rmH_{-1}}{=} \frac{1}{2} \bbE_{\rmP} \left[ \int_r^T|\Psi_s(X_s)|^2 \De s\right] -\bbE_{\rmP}\left[\int_{r}^T \big(\nabla \log \rmP_{s}(X_s)+2\nabla W\ast  \rmP_{s}(X_s) \big) \cdot v_s(X_s)\,\De s\right]\\  \textstyle 
\stackrel{\eqref{chainrule}}{=}\frac{1}{2}\bbE_{\rmP}\left[ \int_{r}^{T} |\Psi_s(X_s)|^2\De s\right]+\nent(\rmP_r)-\nent(\mu^{\mathrm{fin}}).
\eea
The bound \eqref{revcorrectionbound} follows by plugging in \eqref{revcorrectionboundlemma4} into  \eqref{revcorrectionboundlemma2} using the above equation, \eqref{entropytimereversal1} and recalling that $\cH(\rmP|\Gamma(\rmP))=\nonlinentcost(\mui,\muf)$. The proof of \eqref{revcorrectionbound3} goes along the same lines: Since $\hat{\rmP}$ is optimal for \eqref{timerevnonlinearSP} we also get from Lemma \ref{correctionboundlemma}, and in particular from \eqref{correctionbound3} for the choice $t=T-r$ that 
\bes \textstyle 
\frac{1}{2}\bbE_{\hat{\rmP}}\left[ |\hat{\Psi}_{T-r}(X_{T-r})|^2 \right] \leq  \frac{2\curv\, \nonlinentcost(\mu^{\mathrm{fin}},\mu^{\mathrm{in}})}{\exp(2\curv r)-1}.
\ees
\end{proof}


Now the proof of Theorem \ref{nonlinearentropybound1intro} and its corollaries in the introduction is an easy task, given all the preparatory work. 

\begin{proof}[Proof of Theorem \ref{nonlinearentropybound1intro}]
It amounts to add \eqref{correctionbound1} and \eqref{revcorrectionbound} with the choice $r=t$, and use the relation
\bes  \textstyle 
\cH(\rmP|\Gamma(\rmP))=\frac{1}{2}\bbE_{\rmP}\left[\int_0^T|\Psi_t(X_t)|^2\De t\right ] = \nonlinentcost(\mu^{\mathrm{in}},\mu^{\mathrm{fin}}).
\ees
\end{proof}

\begin{proof}[Proof of Corollary \ref{talagrand}]
It follows from Theorem \ref{nonlinearentropybound1intro} (Eq.\ \eqref{nonlinearentropyboundintro}), observing that $\nent(\rmP_t)\geq 0 $.
\end{proof}

\begin{proof}[Proof of Corollary \ref{HWI}]
Combining \eqref{timerevcorrection5},\eqref{timerevcorrection6},\eqref{consquantity2intro} we get that 
\bes \textstyle 
\int_{\RD} |v_t|^2(x) \rmP_t(\De x) = - \cq(\mui,\mu_{\infty}) + \frac{1}{4} \mathcal{I}_{\nent}(\rmP_t).
\ees
Using the above relation, Cauchy Schwartz inequality and the continuity of $\mathcal{I}_{\nent}(\rmP_t)$ in a neighborhood of $0$, \eqref{chainrule} we get that
\be\label{HWI proof eq 1} \textstyle 
\liminf_{t\rightarrow 0} \frac{1}{t}(\nent(\rmP_t)-\nent(\rmP_0) ) \geq - \left( \nonlinfish(\mui) \big( \frac{1}{4}\nonlinfish(\mui)- \cq(\mui,\mu_{\infty})  \big) \right)^{1/2}. 
\ee
Consider now the bound \eqref{nonlinearentropyboundintro}. Observing that $\nent(\mu_{\infty})=0$, subtracting $\nent(\mui)$ on both sides, dividing by $t$, letting $t\rightarrow 0$, using \eqref{HWI proof eq 1} and finally rearranging the resulting terms we get \eqref{HWI consquant}.
\end{proof}

\subsubsection{Proof of Theorem \ref{consquantity1intro}}
We prove here Theorem \ref{consquantity1intro} of the introduction. 
In the proof we will write
\[  \textstyle \int_{\RD} \nabla^2 W(\hat{X}_t-y)\cdot(\hat{\Psi}_t(\hat{X}_t) - \hat{\Psi}_t(y) ) \, \hat{\rmP}_t(\De y),\]
instead of 
\[ \textstyle \tilde{\bbE}_{\tilde{\rmP}} \left[ \nabla^2 W(X_s-\tilde{X}_s)\cdot(\Psi_s(X_s) - \Psi_s(\tilde{X}_s) )\right],\]
which is used in the rest of the article. This is done in order to better deal with time reversal.

\begin{proof}[Proof of Theorem \ref{consquantity1intro}]
Let $M_t$ be the martingale defined at \eqref{martingalecorrection1intro}. Since $\hat{\rmP}$ is optimal for \eqref{timerevnonlinearSP}, from Proposition \ref{martingalecorrection} we get that 

\bes  \textstyle 
\hat{M}_t=\hat{\Psi}_t(X_t)- \int_{0}^t \int_{\RD} \nabla^2 W(\hat{X}_s-y)\cdot(\hat{\Psi}_s(X_s) - \hat{\Psi}_s(y) ) \, \hat{\rmP}_s(\De y) \, \De s
\ees
is an $L^2$-martingale on $[0,T[$ under $\hat{\rmP}$. We define the stochastic processes
\bes  \textstyle 
A_t := \int_{\RD} \nabla^2 W(X_s-y)\cdot(\Psi_s(X_s) - \Psi_s(y) ) \, \rmP_s(\De y) 
\ees
and 
\bes \textstyle 
\hat{A}_t:=\int_{\RD} \nabla^2 W(\hat{X}_s-y)\cdot(\hat{\Psi}_s(\hat{X}_s) - \hat{\Psi}_s(y) ) \, \hat{\rmP}_s(\De y). 
\ees
We have, using the Markovianity of both $\rmP$ and $\hat\rmP$, that
\beas  \textstyle 
\bbE_{\rmP}\left[ \Psi_t(X_t)\cdot \hat{\Psi}_{T-t}(\hat{X}_{T-t}) \right] &=& \textstyle  \bbE_{\rmP}\left[ (M_t +\int_{0}^t A_s \De s) \cdot (\hat{M}_{T-t} +\int_{0}^{T-t} \hat{A}_{s} \De s ) \right]\\
&=&  \textstyle \bbE_{\rmP}\left[ \left(\bbE_{\rmP}[M_T | X_{[0,t]}] +\int_{0}^t A_s \De s\right) \cdot \left(\bbE_{\rmP}[\hat{M}_{T}|\hat{X}_{[0,T-t]} ] +\int_{0}^{T-t} \hat{A}_{s} \De s \right) \right]\\
&=& \textstyle \bbE_{\rmP}\left[ \bbE_{\rmP}[\Psi_T(X_T)-\int_{t}^TA_s\De s\, | X_{[0,t]}] \cdot \bbE_{\rmP}[\hat{\Psi}_T(\hat{X}_T) -\int_{T-t}^{T} \hat{A}_{s} \De s\, |\hat{X}_{[0,T-t]} ] \right]\\
&=&  \textstyle \bbE_{\rmP}\left[ \bbE_{\rmP}[\Psi_T(X_T)-\int_{t}^TA_s\De s\, | X_t] \cdot \bbE_{\rmP}[\hat{\Psi}_T(\hat{X}_T) -\int_{T-t}^{T} \hat{A}_{s} \De s\, |\hat{X}_t ]\right] \\
&=& \textstyle \bbE_{\rmP}\left[\left(\Psi_T(X_T)-\int_{t}^T A_s\De s\right) \cdot \left(\hat{\Psi}_T(\hat{X}_T) -\int_{T-t}^{T} \hat{A}_{s} \De s\,\right) \right].
\eeas
Therefore,
\beas  \textstyle 
\frac{\De}{\De t}\bbE_{\rmP}\left[ \Psi_t(X_t)\cdot \hat{\Psi}_{T-t}(\hat{X}_{T-t}) \right] = \bbE_{\rmP}\left[-A_t \cdot(  \hat{\Psi}_T(\hat{X}_T)-\int_{T-t}^{T} \hat{A}_{s}) + \hat{A}_{T-t}\cdot (\Psi_T(X_T)-\int_{t}^{T} A_s \De s)\right]\\  \textstyle 
=\bbE_{\rmP}\left[-A_t \cdot( \hat{M}_T-\hat{M}_{T-t} + \hat{\Psi}_{T-t}(\hat{X}_{T-t})) + \hat{A}_{T-t}\cdot (M_T-M_t+ \Psi_t(X_t) ) \right].
\eeas
Taking conditional expectation w.r.t.\ $\sigma(X_{[0,t]})$ and using that both $A_t$ and $\hat{A}_{T-t}$ are $X_{[0,t]}$-measurable we get that the above expression equals
\bes \textstyle 
\bbE_{\rmP}\left[-A_t \cdot \hat{\Psi}_{T-t}(\hat{X}_{T-t}) + \hat{A}_{T-t}\cdot \Psi_t(X_t)  \right].
\ees
Using the fact that $W$ is symmetric and the definition of $A_t,\hat{A}_{T-t}$, one easily obtains that the latter expression is worth $0$. Indeed it holds that 
\beas  \textstyle  \bbE_{\rmP}\left[A_t \cdot \hat{\Psi}_{T-t}(\hat{X}_{T-t})\right]=\bbE_{\rmP}\left[\hat{A}_{T-t}\cdot \Psi_t(X_t)  \right]=\\  \textstyle \int_{\RD\times\RD} (\hat{\Psi}_{T-t}(x)-\hat{\Psi}_{T-t}(y))\cdot \nabla^2W(x-y)\cdot (\Psi_{t}(x)-\Psi_{t}(y)) \rmP_t(\De x)\rmP_t(\De y).
\eeas
The proof that the function \eqref{consquantity2intro} is constant on $(0,T)$ is now concluded. In order to establish \eqref{talagrandconsqtyintro} we set $t=T/2$ in \eqref{consquantity2intro} and Cauchy Schwartz to get that \bes  \textstyle  |\cq(\mu^{\mathrm{in}},\mu^{\mathrm{fin}})| \leq \left(\bbE_{\rmP}[|\Psi_{T/2}(X_{T/2})|^2]\bbE_{\rmP}[|\hat{\Psi}_{T/2}(\hat{X}_{T/2})|^2]\right)^{1/2}.\ees The desired conclusion follows from \eqref{correctionbound3} and \eqref{revcorrectionbound3}.
\end{proof}

\subsection{Convergence to equilibrium: proofs}\label{con eq proofs}
\subsubsection{Proof of Theorem \ref{exp conv eq}}

\begin{proof}[Proof of Theorem \ref{exp conv eq}]
Lemma \ref{OTlemma} provides with 
\begin{align*}  \textstyle 
\frac{\De}{\De t} \nent(\rmP_t) & \textstyle \stackrel{\eqref{chainrule}}{=} \bbE_{\rmP}\left[\left( \nabla \log \rmP_t(X_t) + 2\nabla W\ast\rmP_t(X_t)\right)\cdot v_t(X_t)\right] \\ & \textstyle 
\stackrel{\eqref{timerevcorrection5}+\eqref{timerevcorrection6}}{=} \frac{1}{2} \bbE_{\rmP}\left[\left(\Psi_t(X_t)+\hat{\Psi}_{T-t}(\hat{X}_{T-t}) \right)\cdot \left( \Psi_t(X_t)-\hat{\Psi}_{T-t}(\hat{X}_{T-t})\right)\right]\\ & \textstyle 
 =\frac{1}{2}\bbE_{\rmP}\left[|\Psi_t(X_t)|^2-|\hat{\Psi}_{T-t}(\hat{X}_{T-t})|^2\right].
\end{align*}
Going along the same lines as the proof of Lemma \ref{correctionboundlemma} we get that both $\bbE_{\rmP}\left[|\Psi(t,X)|^2 \right]$ and $\bbE_{\rmP}\left[|\hat{\Psi}_{T-t}(\hat{X}_{T-t})|^2 \right]$ are differentiable as functions of $t$ in the open interval $(0,T)$. Moreover 
\beas \textstyle 
\frac{1}{2}\frac{\De}{\De t}\bbE_{\rmP}\left[|\Psi_t(X_t)|^2-|\hat{\Psi}_{T-t}(\hat{X}_{T-t})|^2\right] \stackrel{\eqref{correctionbound5}+\eqref{curvaturebound}}{\geq}  \curv \bbE_{\rmP}\left[|\Psi_t(X_t)|^2+|\hat{\Psi}_{T-t}(\hat{X}_{T-t})|^2\right]\\  \textstyle 
=\curv\, \bbE_{\rmP}\left[|\Psi_t(X_t)+\hat{\Psi}_{T-t}(\hat{X}_{T-t})|^2 -2\Psi_t(X_t)\cdot\hat{\Psi}_{T-t}(\hat{X}_{T-t})\right]\\ \textstyle 
\stackrel{\eqref{timerevcorrection4}}{=}\curv \,\bbE_{\rmP}\left[ |\nabla \log \rmP_t(X_t) + 2\nabla W\ast\rmP_t(X_t)|^2 \right] -2\curv\, \cq(\mu^{in},\mu^{fin}).
\eeas
The $\curv$-convexity of $W$ and the fact that the center of mass $\bbE_{\rmP}[X_t]$ is constant (see Lemma \ref{linearbar}) allow to use the logarithmic Sobolev inequality \cite[(ii), Thm 2.2]{carrillo2003kinetic} to obtain\footnote{ Again, the apparent mismatch between the constant in the Log Sobolev inequality from \cite{carrillo2003kinetic} and the one we use here is due to the fact that in our definition of $\tilde{\cF}$, there is no $1/2$ in front of $W$.}
\bes
\curv \bbE_{\rmP}\left[ |\nabla \log \rmP_t(X_t) + 2\nabla W\ast\rmP_t(X_t)|^2 \right] \geq 4\curv^2 \nent(\rmP_t) .
\ees
Thus we have proven that for almost every $t\in[0,T]$
\bes \textstyle 
\frac{\De t}{\De t^2} \nent(\rmP_t)\geq 4\curv^2\nent(\rmP_t) -2\curv\cq(\mu^{in},\mu^{fin}),
\ees
hence
\eqref{consquantconvexity1intro} follows by integrating this differential inequality (see Lemma \ref{differentialinequality2}). Setting $t=\theta T$ in \eqref{consquantconvexity1intro} and using \eqref{talagrandconsqtyintro} to bound the conserved quantity gives \eqref{exp entropy bound} after some standard calculations.
\end{proof}

\subsubsection{Proof of Theorem \ref{thm exp mkv intro}}

\begin{proof}[Proof of Theorem \ref{thm exp mkv intro}]
Let $\rmP$ be optimal for \eqref{nonlinearSP} and $\Psi$ be given by Proposition \ref{Markovopt}. Then if we define 
\bes \textstyle 
B_t:=X_t - \int_0^t \nabla W \ast \rmP_s(X_s) + \Psi_s(X_s)\De s 
\ees
the process $(B_t)_{t\in[0,T]}$ is a Brownian motion under $\rmP$. Since the McKean Vlasov SDE admits a unique strong solution, there exists a map $\mathbb{Y}:\Omega\longrightarrow\Omega$ such that $\mathbb{Y}\circ B_{\cdot}:=Y$ satisfies 
$Y_0=X_0\,(\rmP-\text{a.s.}) $ and
\bes  \textstyle 
\rmP-\text{a.s.} \quad Y_t = Y_0 - \int_{0}^t \nabla W \ast \rmP^{\text{\tiny{MKV}}}_s(Y_s) \De s + B_t.
\ees
Define now $\delta(t)=\bbE_{\rmP}[|X_t-Y_t|^2]$. Using It\^o's formula we get that $\delta(t)$ is differentiable with derivative
\bes
\delta'(t) = -2\bbE_{\rmP}[(X_t-Y_t)\cdot (\nabla W\ast \rmP_t(X_t)-\nabla W\ast \rmP^{\text{\tiny{MKV}}}_t(Y_t) )] + 2\bbE_{\rmP}[(X_t-Y_t)\cdot \Psi_t(X_t) ] 
\ees

The same arguments as in Lemma \ref{correctionboundlemma} give
\bes
2\bbE_{\rmP}[(X_t-Y_t)\cdot (\nabla W\ast \rmP_t(X_t)-\nabla W\ast \rmP^{\text{\tiny{MKV}}}_t(Y_t) )] \geq 2\curv \bbE_{\rmP}[|X_t-Y_t|^2] \geq 0.
\ees
Moreover, by Cauchy Schwartz:
\bes
 \bbE_{\rmP}[(X_t-Y_t)\cdot\Psi_t(X_t)]  \leq \bbE_{\rmP}\left[|X_t-Y_t|^2\right]^{1/2}\bbE_{\rmP}\left[|\Psi_t(X_t)|^2\right]^{1/2}.
\ees
Therefore
\bes
\delta'(t) \leq 2 \delta(t)^{1/2}\bbE_{\rmP}\left[|\Psi_t(X_t)|^2\right]^{1/2}
\ees
which gives
\bes
(\sqrt{\delta})'(t) \leq  \bbE_{\rmP}\left[|\Psi_t(X_t)|^2\right]^{1/2}.
\ees
Integrating the differential inequality and using that $\delta(0)=0$:
\beas\notag  \textstyle 
\sqrt{\delta}(t) \leq  \int_{0}^t \bbE_{\rmP}\left[|\Psi_s(X_s)|^2\right]^{1/2} \De s \leq t^{1/2} \Big( \int_{0}^t \bbE_{\rmP}\left[|\Psi_s(X_s)|^2\right]\De s\Big)^{1/2} \\ \textstyle 
\stackrel{\eqref{correctionbound1}}{\leq}  t^{1/2} \Big( 2
\frac{\exp(2\curv t)-1}{\exp(2\curv T)-1} \, \nonlinentcost(\mu^{\mathrm{in}},\mu^{\mathrm{fin}}) 
\Big)^{1/2} .
\eeas
The conclusion follows  from \eqref{talagrand1intro} and the observation that $\mathcal W_2^2(\rmP_t,\rmP^{\text{\tiny{MKV}}}_t) \leq \delta(t)$.
\end{proof}

\section{Appendix}

The following lemma is well known. For simple proofs see \cite{LeoGir} or the appendix of \cite{fischer2014form}.

\begin{lemma}\label{wienerent}
For $\mu\in \cP_2(\RD)$ let $\rmR^\mu$ be the law of the Brownian with initial law $\mu$. For $\rmP \in \cP(\Omega)$ with $\rmP\circ(X_0)^{-1}= \mu$ the following are equivalent
\bei 
\item[(i)] $\cH(\rmP|\rmR^\mu)<+\infty$.

\item[(ii)] There exist a $\rmP$-a.s. defined adapted process $(\alpha_t)_{t\in[0,T]}$ such that
\bes\label{finiteent6} \textstyle 
\bbE_{\rmP} \left[ \int_{0}^T |\alpha_t|^2 \De t \right] <+\infty 
\ees
and 
\bes \textstyle 
X_t - \int_{0}^t \alpha_s\De s
\ees
is a Brownian motion under $\rmP$.
\eei
Moreover, if (ii)  holds, then 
\be\label{finitent7} \textstyle 
\cH(\rmP | \rmR^\mu) = \frac{1}{2}\bbE_{\rmP} \left[ \int_{0}^T|\alpha_t|^2 \De t \right].
\ee
and 
\be\label{secondmoment2} \textstyle 
\bbE_{\rmP}\Big[\sup_{t\in[0,T]}|X_t|^2\Big]<+\infty.
\ee
\end{lemma}


With the help of Lemma \ref{wienerent} we can readily prove its generalization given in Lemma \ref{finitent} concerning the case when $R^\mu$ is replaced by the law of a diffusion with Lipschitz drift.

\begin{proof}[Proof of Lemma \ref{finitent}]
The proof of $(i) \Rightarrow (ii)$ follows the arguments in \cite{LeoGir}. 
Now assume that $(ii)$ holds. Because of the continuity of $t\mapsto \bar  b(t,0)$ and \eqref{lipschitz} we get that
\be\label{lineargrowth}
\forall (t,x)\in [0,T] \times \RD, \quad |\bar b(t,x)|\leq C'(1+|x|)
\ee
for some $C'<+\infty$. Consider the sequence of stopping times $T_{n}=\inf\{ t \geq 0 : |X_{t}|=n\}\wedge T$. Using \eqref{lineargrowth},\eqref{finitent12} and some standard calculations we find that there exist $C''<+\infty$ such $\rmP$ almost surely
\bes 
\forall n\in \N, t\in[0,T] :\,\sup_{r\leq t} |X_{r\wedge T_n}|^2 \leq  C''\left(|X_0|^2+1+\int_{0}^{t}\sup_{r\leq s}|X_{r\wedge T_n}|^2 \De s +\int_0^t|\bar\alpha_s|^2 \De s + \sup_{r\leq t}|B_{r\wedge T_n}|^2 \right)
\ees
where $B$ is a Brownian motion. Taking expectation, using \eqref{finitent2}, using Gr\"onwall's Lemma, and eventually letting $n\rightarrow+\infty$, one gets that \be\label{finitent10} \textstyle 
\sup_{t\in[0,T]} \bbE_{\rmP}\Big[\sup_{t\in[0,T]}|X_t|^2\Big] <+\infty.
\ee 
 But then, thanks to \eqref{lineargrowth} and \eqref{finitent9} we also obtain that 
$\bbE_{\rmP}[\int_{0}^T |\bar b(s,X_s) + \bar\alpha_s|^2\De s  ]<+\infty$. Lemma \ref{wienerent} yields  then (as usual $\rmR^\mu$ is Wiener measure started like $\mu$)
\be\label{finitent4} \textstyle 
\cH(\rmP|\rmR^\mu)<+\infty, \quad \text{and} \quad\cH(\rmP|\rmR^\mu) = \frac{1}{2}\bbE_{\rmP}\left [\int_{0}^T |\bar b(t,X_t) + \bar\alpha_t|^2\De t \right ].
\ee
Under the current hypotheses on $b$, $\rmR^\mu$ and $\bar\rmR$ are mutually absolutely continuous and
\bes \textstyle 
\frac{\De \rmR^\mu}{\De\bar \rmR} = \exp\left( -\int_0^T \bar  b(t,X_t) \cdot \De X_t +\frac{1}{2} \int_0^T |\bar b(t,X_t)|^2 \De t \right).
\ees
Therefore, using some standard calculations and $(ii)$ we get
\be\label{finitent3} \textstyle 
{\bbE_{\rmP} \left[ \log \frac{\De \rmR^\mu}{\De\bar \rmR} \right] }= -\bbE_{\rmP}\left[\int_{0}^T \left(\bar\alpha_t+ \frac{1}{2}\bar b(t,X_t)\right)\cdot \bar b(t,X_t) \De t \right]<+\infty.
\ee
Since $\rmR^\mu$ and $\bar \rmR$ are mutually absolutely continuous and $\cH(\rmP|\rmR^\mu)<+\infty$ we get that $\rmP \ll\bar \rmR$ and 
\be\label{finitent5} \textstyle 
\cH(\rmP|\bar \rmR) = \bbE_{\rmP}\left[\log \frac{\De \rmP }{\De \rmR^\mu } +\log \frac{\De \rmR^\mu }{\De\bar \rmR } \right]
\ee
Thus, $\cH(\rmP|\bar\rmR)<+\infty$ if both
$ \cH(\rmP|\rmR^\mu)$ and $\bbE_{\rmP}\left[\log \frac{\De \rmR^\mu}{\De \bar\rmR}\right]$ are finite. But this is indeed the case, thanks to \eqref{finitent3},\eqref{finitent4}. The proof that $(ii)\Rightarrow(i)$ is now complete.  The desired form of the relative entropy follows by plugging in \eqref{finitent4} and \eqref{finitent3} into \eqref{finitent5}.  Finally \eqref{secondmoment} follows from \eqref{finitent10} and \eqref{lineargrowth}, and the last statement from \eqref{secondmoment}.
\end{proof}

\begin{lemma}\label{entropyequivalence}
Let $\mu^{in}\in\cP(\RD)$ and $\rmP\in\cP_{1}(\Omega)$ with $\rmP_0=\mui$. Then $\cH(\rmP | \Gamma(\rmP))<+\infty$ if and only if $\cH(\rmP |  \rmR^{\mui} )<+\infty$.
\end{lemma}

\begin{proof}
 Define $\bar{b}(t,z) = -\nabla W \ast \rmP_t(z)$. Lemma \ref{driftregularity} ensures that $\bar{b}$ is of class $\cC^{0,1}$ and that \eqref{lipschitz} holds. Assume that $\cH(\rmP | \Gamma(\rmP))<+\infty$ and let $(\bar{\alpha}_t)_{t\in[0,T]}$ be given by the implication $(i)\Rightarrow (ii)$ of Lemma \ref{finitent}. If we define $(\alpha_t)_{t\in[0,T]}$ by
  \bes \rmP-a.s. \quad \alpha_t = \bar{b}(t,X_t)+ \bar{\alpha}_t, \quad \forall t\in[0,T], \ees
  then \eqref{secondmoment} together with \eqref{finitent9} entitle us to use the implication $(ii)\Rightarrow(i)$ of Lemma \ref{wienerent} to obtain the desired result. The converse implication is done inverting the roles of Lemmas \ref{finitent} and \ref{wienerent}.
\end{proof}

\begin{lemma}\label{entropyapprox}
Assume \eqref{boundedhess},\eqref{marginalhyp} $\rmP,\rmQ\in\cP_{1}(\Omega)$ be such that $\cH(\rmP|\Gamma(\rmP)),\cH(\rmQ|\Gamma(\rmP))<+\infty$ and $(\rmQ_t)_{t\in[0,T]}=(\rmP_t)_{t\in[0,T]}$. If $(\theta^n)_{n\in\N} \subseteq \cC^{\infty}_{c}([0,T]\times\RD)$  is such that $\nabla \theta^n$ converges to $-\nabla W \ast \rmP_t(z)$ in $\rmH_{-1}((\rmP_t)_{t\in[0,T]})$, i.e.
\be\label{L2conv} \textstyle 
\lim_{n\rightarrow +\infty}  \int_{[0,T]\times \RD} |\nabla \theta^n_t(z)+\nabla W \ast \rmP_t(z)|^2  \rmP_t(\De z)\,\De t  =0,
\ee
then 
\bes
\lim_{n\rightarrow +\infty} \cH(\rmQ|\rmR^n)= \cH(\rmQ|\Gamma(\rmP)),
\ees
where $\rmR^n$ is the law of
\bes
\De Y_t = \nabla \theta^n_t(Y_t)\De t + \De B_t, \quad Y_0 \sim \mui\in\cP_2(\mathbb R^d).
\ees
\end{lemma}

\begin{proof}
By Lemma \ref{entropyequivalence} we have $\cH(\rmQ|\rmR^{\mui})<+\infty$, where $R^{\mui}$ is the law of the Brownian motion with initial law $\mui$. Using implication $(i)\Rightarrow(ii)$ from Lemma \ref{wienerent} and then $(ii)\Rightarrow (i)$ together with \eqref{finitent13} from Lemma \ref{finitent} for the choice $\bar b=\nabla \theta^n$, we get for all $n\in \N$:
\be\label{entopyapprox4}
\cH(\rmQ|\rmR^n) = \frac{1}{2}\bbE_{\rmQ} \left[\int_0^T |\alpha_t - \nabla \theta^n_t(X_t)|^2 \De t \right],
\ee
where $\alpha_t$ is the dirft of $\rmQ$ (see Lemma \ref{wienerent}).
Moreover, using $\cH(\rmQ|\Gamma(\rmP))<+\infty$ and Lemma \ref{finitent}, we also get that 
\be\label{entopyapprox5} \textstyle 
\cH(\rmQ|\Gamma(\rmP)) = \frac{1}{2}\bbE_{\rmQ} \left[\int_0^T |\alpha_t + \nabla W\ast \rmP_t(X_t) |^2 \De t \right].
\ee
Using $(\rmQ_t)_{t\in[0,T]}=(\rmP_t)_{t\in[0,T]}$ and \eqref{L2conv} we get 
\be\label{entopyapprox1} \textstyle 
\lim_{n\rightarrow+\infty}\bbE_{\rmQ} \left[\int_0^T |\nabla \theta^n_t(X_t)|^2 \De t \right] = \bbE_{\rmQ} \left[\int_0^T |\nabla W\ast \rmP_t(X_t)|^2 \De t \right].
\ee
On the other hand, let $\bar{\alpha}_t(X_t)$ be a measurable version of $\bbE_{\rmQ}[\alpha_t| X_t] $, the existence of which is guaranteed e.g.\ by \cite[Proposition 4.4]{gyongy1986mimicking}.
Using conditional Jensen and $\bbE_{\rmQ}\left[ \int_{0}^T|\alpha_t|^2\De t\right]<+\infty$, it is easy to verify that $\bar{\alpha} \in \rmH_{-1}((\rmP_t)_{t\in[0,T]})$. Moreover
\bes \textstyle 
\bbE_{\rmQ} \left[\int_0^T\alpha_t \cdot \nabla \theta^n_t(X_t) \De t \right] = \int_{[0,T] \times \RD} \bar{\alpha}_t(z)\cdot \nabla \theta^n_t(z)  \rmP_t(\De z) \De t.
\ees
Since $\bar{\alpha} \in \rmH_{-1}((\rmP_t)_{t\in[0,T]})$ we get from \eqref{L2conv}, $(\rmQ_t)_{t\in[0,T]}=(\rmP_t)_{t\in[0,T]}$ and the basic properties of conditional expectation 
\be\label{entopyapprox2} \textstyle 
\lim_{n\rightarrow +\infty} \int_{[0,T] \times \RD} \bar{\alpha}_t(z)\cdot \nabla \theta^n_t(z)  \rmP_t(\De z) \De t=\bbE_{\rmQ}\left[ \int_{0}^T \alpha_t \cdot \nabla W \ast \rmP_t(X_t)\De t\right].
\ee
%
Gathering \eqref{entopyapprox5},\eqref{entopyapprox2},
\eqref{entopyapprox4},\eqref{entopyapprox1} the conclusion follows.
\end{proof}

\begin{lemma}\label{shiftinvert}
Assume 
$\rmP\ll\rmR^\mu$ and that $(h_t)_{t\in[0,T]}$ is an adapted process satisfying 
the Lipschitz condition in \eqref{Lipschitzperturb}. Then the shift $\tau_{\varepsilon}$ defined at \eqref{shiftdef} admits an almost sure inverse, i.e.\ there exists an adapted process $(Y^{\varepsilon}_{t})_{t\in[0,T]}$ such that 
\be\label{shiftinvert1}
 \rmP-\text{a.s.}\quad \tau^{\varepsilon}_t (Y^{\varepsilon} (\omega))=Y^{\varepsilon}_t (\tau^{\varepsilon} (\omega)) = \omega_t \quad \forall t\in [0,T].
\ee
\end{lemma}

\begin{proof}
Let $\rmR^\mu$ be the law of the Brownian motion started at $\mu$. The fact that  \eqref{shiftinvert1} holds $\rmR^\mu$ almost surely under the Lipschitz condition \eqref{Lipschitzperturb} is a classical result, see e.g. \cite[pg 209-210]{ustunel2007sufficient}. In this case the a.s.\ inverse is nothing but the strong solution of the SDE 
\bes
\De Y^{\varepsilon}_t = -\varepsilon h_t(Y^{\varepsilon}_t)\De t + \De B_t
\ees
We conclude by $\rmP\ll\rmR^\mu$.
\end{proof}

The next lemma follows from \cite[Lemma 4.1]{conforti2018second}:

\begin{lemma}\label{differentialinequality}
	Let $c:[0,T]\rightarrow \R$ be twice differentiable on $(0,T)$ with $c(0)=0$ and let $\curv\in\R$.  If $\frac{\De^2}{\De t^2}c(t) \geq 2\curv \frac{\De}{\De t}c (t)$ for all $t \in (0,T)$, then
	 \be \textstyle \label{eq:diffineq1} \forall t \in [0,T], \quad c(t) \leq  \frac{\exp(2\curv t)-1}{\exp(2\curv T)-1} c(T).\ee
\end{lemma}

The following lemma is also useful for the quantitative estimates:

\begin{lemma}\label{differentialinequality2}
	Let $\varphi:[0,T]\rightarrow \R$ be twice differentiable on $(0,T)$ and assume that $\varphi(t)$ satisfy the differential inequality
\bes \textstyle  
\frac{\De^2}{\De t^2} \varphi(t) \geq \lambda^2 \varphi(t)-\lambda \, \mathscr{E}
\ees
where $\mathscr{E}$ is a constant. Then we have for all $0\le t\le T,$ 
\be\label{comparison3} \textstyle 
\varphi(t) \leq (\varphi(0)-\frac{\mathscr{E}}{\lambda}) \alpha_{\lambda,T}(T-t) + (\varphi(T)-\frac{\mathscr{E}}{\lambda})\alpha_{\lambda,T}(t) + \frac{\mathscr{E}}{\lambda}. \ee
 
\end{lemma}

\begin{proof}[Proof of Lemma \ref{differentialinequality2}]
Consider the function $\gamma_t= \varphi_t-\psi_t$, where $\psi_t$ is the right hand side of \eqref{comparison3}. Observe that $ \frac{\De^2}{\De t^2}\gamma_t \geq \lambda^2 \gamma_t$ for $t\in[0,T]$ and $\gamma_0=\gamma_T=0$.  Assume $\gamma_{t_1}>0$ for some $t_1$. Since $\gamma_{0}=0$, there must exist $t_0\leq t_1$ such that $\gamma_{t_0}>0$ and $\frac{\De}{\De t} \gamma_{t_0}>0$. But this contradicts Lemma \ref{comparison2} below.
\end{proof}

\begin{lemma}\label{comparison2}
Assume that $ \lambda\ge 0.$
Let $0\leq t_0< T$  and $\gamma:[0,T]\rightarrow \R$ be a function satisfying 
$$ \bec \frac{\De^2}{\De t^2} \gamma_t \geq \lambda \gamma_t,   \mbox{ $t\in [t_0,T],$} \\ \gamma_{t_0}>0, \gamma_T = 0. \eec $$
Then $\frac{\De}{\De t}{\gamma}_{t_0} \leq 0$.
\end{lemma}

\begin{proof} 
Assume ad absurdum  that $\frac{\De}{\De t}{\gamma}_{t_0} > 0$. Since $\gamma_T=0$ it must be that $t_1=\inf\{t\geq t_0 ; \frac{\De}{\De t} \gamma_t =0 \}$ belongs to $(t_0,T]$. By definition of $t_1$, and since $\frac{\De}{\De t}{\gamma}_{t_0} > 0$ we have that $\gamma_t\geq \gamma_{t_0}>0$ for all $t\in [t_0,t_1]$. But then 
\beas \textstyle 
0> \frac{\De}{\De t}{\gamma}_{t_1}-\frac{\De}{\De t}{\gamma}_{t_0} = \int_{t_0}^{t_1} \frac{\De^2}{\De t^2}{\gamma}_{t} \De t \geq  \lambda \int_{t_0}^{t_1} {\gamma}_{t} \De t\ge 0,
\eeas
which is absurd.
\end{proof}

\bibliographystyle{plain}
\bibliography{Ref}
\end{document}